%% file: time_freezing_final.tex
\documentclass[letterpaper, 10 pt, conference]{ieeeconf}  % Comment this line out
\IEEEoverridecommandlockouts                              % This command is only

\usepackage{amsmath,amssymb,mathtools,amsfonts}
\usepackage[pdftex]{graphicx}
\usepackage{enumerate}
\usepackage{dsfont}
\usepackage{comment}

\usepackage[normalem]{ulem}

\usepackage{multirow}
\usepackage{graphicx}      % include this line if your 
\usepackage{color}
\usepackage{psfrag}
\usepackage{subfigure}
\usepackage{empheq}
\usepackage{tabularx}
\usepackage{bm} % bold italic symbols in math mode
\usepackage{nicefrac} % for nice inline fractions
\usepackage{url}
\usepackage{psfrag}
\usepackage{multirow}
\usepackage{psfrag}
\usepackage{pgfplots}
\pgfplotsset{compat=newest}
\usepackage{float}
\usepackage{amsthm}
\usepackage{alphalph}
\usepackage{booktabs}

\usepackage[short]{optidef}
\usepackage{tikz}
\usepackage{pgfplots}
\usetikzlibrary{external}
\newlength\fwidth
\newlength\fheight
\definecolor{offwhite}{RGB}{249,242,215}
\definecolor{foreground}{RGB}{0,0,0}
\definecolor{background}{RGB}{255,255,255}
\definecolor{title}{RGB}{0,0,0}
\definecolor{gray}{RGB}{155,155,155}
\definecolor{subtitle}{RGB}{155,155,155}
\definecolor{hilight}{RGB}{102,255,204}
\definecolor{vhilight}{RGB}{255,111,207}
\definecolor{lolight}{RGB}{155,155,155}
\definecolor{blueg}{RGB}{114,143,178}
\definecolor{bluet}{RGB}{100,121,178}
\definecolor{redt}{RGB}{255,86,73}
\definecolor{greent}{RGB}{136,255,134}
\definecolor{dgreent}{RGB}{75,212,69}
\definecolor{lgray}{RGB}{190,190,190}
\definecolor{llgray}{RGB}{220,220,220}

\newif\ifcommentarmin
\commentarminfalse
\commentarmintrue
\newcommand{\armin}[1]{\ifcommentarmin{\textcolor{blueg}{(armin: \textbf{#1})}}\else \fi}

\newtheorem{theorem}{Theorem}

\newtheorem{definition}{Definition}
 
\newtheorem{Assumption}{Assumption}

\newcommand{\R}{{\mathbb{R}}}

\newcommand{\cC}{{\mathcal C}}

\newcommand{\dd}{{\mathrm{d}}}

\title{\LARGE \bf
A Time-Freezing Approach for Numerical Optimal Control of Nonsmooth Differential Equations with State Jumps
}

\author{Armin Nurkanovi\'c, Tommaso Sartor, Sebastian Albrecht, Moritz Diehl%<-this % stops a space
\thanks{This research was supported by the German Federal Ministry of Education and Research (BMBF) via the funded Kopernikus project: SynErgie (03SFK3U0), by the German Federal Ministry for Economic Affairs and Energy (BMWi) via DyConPV (0324166B), and by DFG in project DI 905/3-1 and via Research Unit FOR 2401.}
\thanks{Armin Nurkanovi\'c  and Sebastian Albrecht are with with Siemens Corporate Technology, 81739 Munich, Germany. Armin Nurkanovi\'c  is also with the Department of Microsystems Engineering (IMTEK) University Freiburg, 79110 Freiburg, Germany
{\tt\small \{armin.nurkanovic, sebastian.albrecht\}@siemens.com}}%
\thanks{Tommaso Sartor is with the MECO Research Team, Department Mechanical Engineering, KU Leuven, Leuven, Belgium  {\tt\small \{  tommaso.sartor\}@kuleuven.be}}
\thanks{Moritz Diehl is with the Department of Microsystems Engineering (IMTEK) and Department of Mathematics, University Freiburg,
        79110 Freiburg, Germany
        {\tt\small \{moritz.diehl\}@imtek.uni-freiburg.de}}%
}% <-this % stops a space

\begin{document}

\maketitle
\thispagestyle{empty}
\pagestyle{empty}

%%%%%%%%%%%%%%%%%%%%%%%%%%%%%%%%%%%%%%%%%%%%%%%%%%%%%%%%%%%%%%%%%%%%%%%%%%%%%%%%
\begin{abstract}
We present a novel reformulation of nonsmooth differential equations with state jumps enabling their easier simulation and use in optimal control problems without the need for integer variables. 
The main idea is to introduce an auxiliary differential equation to mimic the state jump map. Thereby, a clock state is introduced which does not evolve during the runtime of the auxiliary system. The pieces of the trajectory that correspond to the parts when the clock state was evolving recover the solution of the original system with jumps.
Our reformulation results in nonsmooth ordinary differential equations where the discontinuity is in the first time derivative of the trajectory, rather than in the trajectory itself. This class of systems is easier to handle both theoretically and numerically. 
We provide numerical examples demonstrating the ease of use of this reformulation in both simulation and optimal control.
\textcolor{black}{In the optimal control example, we solve a sequence of nonlinear programming problems (NLPs) in a homotopy penalization approach and recover a time-optimal trajectory with state jumps.}

\end{abstract}
%Dynamic systems with state jumps can be described by measure differential inclusions (MDIs). Our reformulation gives a specific nonsmooth ordinary differential equations which is easier to handle both theoretically and numerically then an MDI.
%We show how the true initial trajectory can be recovered from the simpler system in a general case.
%%%%%%%%%%%%%%%%%%%%%%%%%%%%%%%%%%%%%%%%%%%%%%%%%%%%%%%%%%%%%%%%%%%%%%%%%%%%%%%%
\section{Problem Description} \label{sec:ProblemDescription}
This paper regards the numerical treatment of nonsmooth differential equations in optimal control. The nonsmoothness of  ${\dot{x}(t) = f(x(t))}$ can be classified depending on the classes into which the solution $x(t;x_0)$ and right hand side (r.h.s.) $f(x(t))$ fall: 1) Ordinary Differential Equations (ODEs) with nonsmooth but Lipschitz r.h.s. and $\cC^1$ solutions; 2) discontinuous but one-sided Lipschitz r.h.s. with absolutely continuous (AC) solutions; 3) solutions that contain state jumps and are functions of bounded variations. This paper focuses on case 3. Since $x(t;x_0)$ jumps, $f(x(t))$ has to contain Dirac-$\delta$ impulses. In such cases we cannot in general speak of ODEs and we have to use tools such as Measure Differential Inclusions (MDIs) \cite{Moreau1977}. 
%\textcolor{red}{More compactly, the question is whether a jump discontinuity appears in $x(t;x_0)$ (case 3) or in its first or second time derivative (cases 2 and 1).}
These differential equations arise in: rigid-bodies with friction and impact, electronics, traffic flows, biological systems, economical systems, energy systems, cf. \cite{Brogliato2020}.

There are many different formalisms to model nonsmooth dynamic phenomena, for an overview the reader is referred to the excellent monographs \cite{Brogliato2020,Acary2008,Stewart2011}. 
Despite the very good and solid developments both in theory (e.g., existence and uniqueness of solutions for various formalisms \cite{Brogliato2020,Stewart2011}) and numerical simulation methods \cite{Acary2008}, there is still a lack of practical numerical optimal control methods for the three mentioned classes of dynamics systems. 
%\textcolor{red}{In this work we focus on the third class above: on nonsmooth differential equations with state jumps.} 
While the first class poses no major obstacle to practical solution, e.g. with smoothing, the second and third classes are difficult.
%
%\textcolor{red}{The solutions of the second and third class are piecewise smooth. A common approach is to use binary indicator variables for the smooth pieces and apply Mixed Integer Nonlinear Programming (MINLP) \cite{Oldenburg2005}. An alternative to MDIs or Differential Variational Inequalities (DVIs) of index 2 \cite{Pang2008} to model dynamics with state jumps is to use the \textit{hybrid systems} formalism \cite{Goebel2009}. Thereby, after the trajectories hit the boundary of a \textit{flow set} an algebraic \textit{jump map} is used for the state reinitialization. OCP formulations with this formalism result also in MINLPs \cite{Altin2018}. Nonconvex MINLP approaches are often far from being computationally tractable.}
\textcolor{black}{To mitigate the difficulties caused by the state jumps, two common approaches are: (a) to use some coordinate transformations \cite{Kim2014,Zhuravlev1978} and (b) to use smoothing/penalization \cite{Tassa2012} or some compliant impact model \cite{Brogliato2016}. The goal is to obtain dynamics which fall into case 2 or even 1. Coordinate transformation can be very efficient for some special settings. The Zhuravlev-Ivanov transformations \cite[Sec. 1.4.3]{Zhuravlev1978,Brogliato2016} are restricted to mechanical systems and to constraints of co-dimension one. A more general approach is the use of \textit{gluing functions} within the \textit{hybrid systems} formalism \cite{Kim2014}. However, this approach regards also only co-dimension one constraints and there is no algorithm for finding the needed gluing function. To obtain realistic approximations with smoothing/penalization one has to deal with very stiff differential equations, and compliant models can yield nonphysical effects \cite[Sec. 2.2]{Brogliato2016}. The method introduced in this paper falls somewhat in between these two approaches, as we also transform the system into an equivalent system which has AC solutions on a different time domain, and since it can use complaint models to emulate state jump laws.\\
Modeling switched systems with complementarity conditions (CC) is gaining more popularity \cite{Baumrucker2009}. Walking, running and manipulation problems are rich sources of  Optimal Control Problems (OCPs) with nonsmooth dynamics in robotics \cite{Posa2014}. Rigid-body impact problems with friction are often modeled via Dynamic Complementarity Systems (DCSs). In few recent papers \cite{Guo2016,Vieira2019} the authors study necessary and sufficient conditions in function spaces for OCPs with DCSs with AC solutions. Discretization of OCPs with CCs results in Mathematical Programs with Complementarity Constraints (MPCCs). Unfortunately, difficulties with numerical sensitivities arise when one discretizes the CCs within direct methods.} Conditions for obtaining the right numerical sensitivities with smoothing of differential equations with a discontinuous r.h.s. are provided in the excellent paper by Stewart and Anitescu \cite{Stewart2010}. Their result is extended to MPCCs originating from OCPs in \cite{Nurkanovic2020}. Many MPCC algorithms use smoothing, relaxation or penalty methods \cite{Ralph2004}. The main conclusion from these papers is: in direct collocation for the case 2 of dynamic systems one has to use a sufficiently small step size in comparison with the smoothing parameter, so that the sensitivities of the smoothed system approach the sensitivities of the nonsmooth dynamic system.

%\subsection{Contributions and Outline}
\paragraph*{\textcolor{black}{Contributions}}In this paper we present a novel formulation of restitution laws for nonsmooth differential equations with state jumps. The main idea is to introduce an auxiliary dynamic system, where the initial and endpoint of the solution on some interval satisfy the restitution law. 
%The time when this dynamic is active we denote as \textit{restitution phase}.
Furthermore, a clock state is introduced which does not evolve when the auxiliary dynamic system is active. Finally, we take the pieces of the trajectories corresponding to the time intervals where the clock state is evolving and thereby we recover the solution of the original dynamic system with state jumps. The efficacy of this approach is demonstrated in both simulation and optimal control experiments.
\paragraph*{\textcolor{black}{Outline}}
The paper is structured as follows: Section \ref{sec:TimeFreezing} introduces the main ideas and all terminology, followed by Section \ref{sec:TimeFreezingExample} where all concepts are illustrated on a simple example. In Section \ref{sec:Equivalence} we relate the solutions of the original dynamic system and our reformulation and show how to recover the original solution in the general case. Section \ref{sec:NumericalExamples} provides both simulation and optimal control examples. \textcolor{black}{We solve a time-optimal control problem of a moving particle where the optimal solution considers multiple simultaneous impacts.}
% We compare our approach to a handcrafted multi-stage formulation and get the best solution without the need of stage number enumeration or heuristics.
 The paper concludes and discusses further extensions in Section \ref{sec:Conclusions}.

\paragraph*{\textcolor{black}{Notation}}
For the time derivative of a function $x(t)$ we use $\dot{x}(t) \coloneqq \frac{dx(t)}{dt} $ and for $y(\tau)$ we use ${y}'(\tau)\coloneqq\frac{dy(\tau)}{d\tau} $. For the left and right limit, we use the notation ${x(t_s^+)  = \lim\limits_{t\to t_s,\ t>t_s}x(t)}$ and ${x(t_s^-)  = \lim\limits_{t\to t_s,\  t<t_s}x(t)}$, respectively. 
%The complementary conditions for two vectors  $a,b \in \R^{n}$ reads as: ${0\leq a \perp b\geq 0}$. 
%The inequalities are to be understood element-wise and $a \perp b$ means $a \circ b = 0$ where $\circ$ is the Hadamard element-wise product (or equivalently $a^{\top}b =0$, since $a,b\geq 0$). }
The matrix ${\mathds{1}_{n,n} \in \R^{n\times n}}$ is the identity matrix, and ${0}_{m,n} \in \R^{m\times n}$ is the zero matrix. \textcolor{black}{A column vector in  $\R^k$ with all ones is denote as $e_k$.
The concatenation of two column vectors $a\in \R^m$, $b\in \R^n$ is denoted as $(a,b)\coloneqq[a^\top,b^\top]^\top$.}
\section{Time-Freezing of Differential Equations with State Jumps} \label{sec:TimeFreezing}
We regard differential equations with unilateral constraints and state jumps: 
$\dot{x}(t) = f(x(t)), \ \psi(x(t)) \geq0 $. The \textit{switching manifold} $S$ is defined as $ S \coloneqq \{x \ |  \ \psi(x) =0 \}$  and splits the state space $\R^{n_x}$ into two pieces: the \textit{feasible region} $V^+ \coloneqq \{x \ | \ \psi(x) \geq 0 \}$ and the \textit{prohibited region } $V^- \coloneqq\{x \ | \ \psi(x) <0 \}$. Moreover, depending in which direction the trajectory points, and using $\dot{\psi}(x)\coloneqq \nabla \psi(x)^{\top} f(x)$, the switching manifold can be split into the following subsets: $S^+ \coloneqq \{x \ | \ x \in S,\ \dot{\psi}(x) >0 \}$, $S^-\coloneqq \{x \ |\ x \in S, \ \dot{\psi}(x)  <0 \}$ and  ${S^0 \coloneqq \{x \ |\ x \in S, \ \dot{\psi}(x)  = 0 \}}$. At time of impact $t_s$, just before the impact $x(t_s^-) \in S^-$ and the trajectory points outside the feasible region $V^+$ ($\nabla \psi(x(t_s^-)^\top f(x(t_s^-)) <0$). To keep the trajectory feasible, a state jump has to occur so that the trajectory points again into $V^+$, i.e. $x(t_s^+) \in S^+$. This is achieved with the restitution law $\Gamma: S^- \to S^+$. We collect these properties in the following definition:
\begin{definition}[Ordinary Differential Equation with State Jumps]\label{def:ODE_state_jump}
	%We define the time $t \in \R$ ,	the differential states $x(t) \in \R^{n_x}$ and control inputs $u(t) \in \R^{n_u} $. A system of differential equations with state jumps describe the dynamic evolution of	this state vector $x(t)$ as
	We define the time $t \in \R$ and the differential states $x(t) \in \R^{n_x}$. A system of differential equations with state jumps describes the dynamic evolution of the state vector $x(t)$ as
\begin{subequations}\label{eq:DifferentialJump}
	\begin{align}
		\dot{x}(t)  &=	f(x(t)),  \ x(t) \in V^+, \label{eq:DifferentialJump_ode}\\ 
		 x(t^+)&= \Gamma(x(t^-)),  \ \textrm{if}\ \psi(x(t)) = 0  \ \textrm{and}\ x(t^-)\in S^- \label{eq:DifferentialJump_jump},
	\end{align}
\end{subequations}
	where  $\psi(x(t)): \R^{n_x} \to \R$ describes a constraint on the dynamics $f:\R^{n_x} \to \R^{n_x}$. The function ${\Gamma : S^- \to S^+}$ is the restitution law and is used at all $t$ where $x(t) \in S^-$.
\end{definition}
%(Given an initial condition $x(t_0) = x_0$, a set of differential equations \eqref{eq:DifferentialJump} is denote as an \textit{initial value problem (IVP)} over a time interval $t\in[t_0,t_f]$.)

In case of mechanical impact problems, such systems are sometimes called \textit{vibro-impact systems} \cite{Brogliato2016}. As an example of such systems we consider the dynamics of a ball bouncing on a table, which is given by:
\begin{subequations}\label{eq:bouncing_ball_model}
\begin{align}
	m{\dot{v}(t)} &= -mg,\	\dot{q}(t) = v(t), \  \textrm{ if }  q(t) \geq 0  \\
	v(t^+) &= -\gamma v(t^-),\ \textrm{whenever}\ q(t) = 0 \land v(t) < 0, \label{eq:restitution_law} 
	\end{align}
\end{subequations}
where $q(t)$ is the height of the ball, $v(t)$ is the velocity of the ball, $m$ is the mass of the ball and $g$ is the gravitational acceleration. Equation \eqref{eq:restitution_law} is Newton's restitution law for impact dynamics, where $\gamma\in\left[0,1\right]$ is the \textit{coefficient of restitution}. Several other restitution laws can be found in the literature, cf. \cite{Brogliato2016}. 
%The goal of the restitution law is to prevent the ball from penetration (i.e. violating $q(t)\geq 0$). 
%Furthermore, it changes the velocity of the ball at time of impact $t$, which is $v(t^-) < 0$ to a value $v(t^+) \geq 0$ instantaneously. 

%We turn back to the general form in Definition \ref{def:ODE_state_jump} and define some more terminology which is used throughout the paper.
%The \textit{switching manifold} $S$ is defined as $ S = \{x(t) \ | \ \psi(x(t)) =0 \}$  and splits the state space $\R^{n_x}$ into two pieces: the \textit{feasible region} $V^+ \coloneqq \{x(t) \ | \ \psi(x(t)) \geq 0 \}$ and a \textit{prohibited region } $V^- \coloneqq\{x(t) \ | \ \psi(x(t)) <0 \}$. Moreover, depending in which direction the trajectory points, the switching manifold can be split into the following subsets: ${S^+ \coloneqq \{x(t) \ | \ \dot{\psi}(x(t)) >0 \}}$, ${S^-\coloneqq \{x(t) \ | \ \dot{\psi}(x(t))  <0 \}}$ and  ${S^0 \coloneqq \{x(t) \ | \ \dot{\psi}(x(t))  = 0 \}}$, where $S = S^+ \cup S^-\cup S^0$. At time of impact $t$, $x(t) \in S^-$ and the trajectory points outside  the feasible region $V^+$ ($\nabla \psi(x)f(x(t) <0$). To keep the trajectory feasible a state jump has to occur so that trajectory points into $V^+$, i.e. $x(t) \in S^+$. This is achieved with the restitution law $\Gamma: S^- \to S^+$. 
 Since in the general case, the time of impact $t$ is not known a priori, simulating and incorporating such models with additional algebraic conditions into optimization problems is difficult. To alleviate all these difficulties we propose the following approach. First, we relax the constraint $\psi(t)\geq 0$ and define an auxiliary dynamic system on $V^-$ to mimic the restitution law. Second, we introduce a clock state $t$ that evolves according to $t'(\tau) = 1$.  The time $\tau$, denoted as \textit{pseudo time} is now the time of the differential equation. The state evolution of $x(\cdot)$ from Definition \ref{def:ODE_state_jump} in pseudo time is denoted as $\tilde{x}(\tau)$. Third, we "freeze" the time whenever $\tilde{x}(\tau)\in V^-$, i.e. ${{t}'(\tau) = 0}$.
% Third, whenever $\tilde{x}(\tau) \in V^+$, the time is allowed to evolve, i.e ${{t}'(\tau) = 1}$. 
To mimic the restitution law, we assume there exists an auxiliary ODE, whose endpoints satisfy the restitution law on a finite time interval:
% $(\tau_{\textrm{s}},\tau_{\textrm{r}})$, i.e. $ \tilde{x}(\tau_{\textrm{r}}) = \Gamma (\tilde{x}(\tau_{\textrm{s}}))$. 
\begin{Assumption}\label{ass:Auxiliary_Dynamics}
There exists an auxiliary dynamic system ${\tilde{x}}'(\tau) = \varphi(\tilde{x}(\tau))$ such that for every initial value ${\tilde{x}(\tau_{\textrm{s}})=\tilde{x}_0 \in S^-}$, the following properties hold on a finite and well-defined time interval $(\tau_{\textrm{s}},\tau_{\textrm{r}}),$ (with ${\tau_{\textrm{jump}} \coloneqq \tau_{\textrm{r}} -\tau_{\textrm{s}}}$)$\ \colon$
$\tilde{x}(\tau) \in V^-,\ \forall \tau\in(\tau_{\textrm{s}},\tau_{\textrm{r}})$, the dynamics \textcolor{black}{has} its first intersection with $S$ after $\tau_{\textrm{jump}}$ with $\tilde{x}(\tau_{\textrm{r}}) \in S^+$ and ${\tilde{x}(\tau_{\textrm{r}}) = \Gamma (\tilde{x}(\tau_{\textrm{s}}))}$.
\end{Assumption}
%We collect all the introduced ideas in the following definition:
\noindent\textcolor{black}{The introduced ideas are collected in the following definition.}
\begin{definition}[Time-Frozen Differential Equations]\label{def:ODE_time_frozeon}
We define the pseudo-time $\tau \in \R$, the differential states ${y(\tau)\coloneqq (\tilde{x}(\tau),t(\tau)) \in \R^{n_x+1}}$. A system of differential equations  describes the dynamic evolution of the state vector $y(\tau)$ as
%\begin{subequations}\label{eq:DifferentialTimeForezen}
%	\begin{align}
%		{y}'(\tau) &=\begin{cases}
%		\tilde{f}(y(\tau)),  \ &\tilde{\psi}(y(\tau)) \geq 0,\label{eq:DifferentialTimeForezen_ode} \\
%		\tilde{\varphi}(y(\tau)) ,  \ &\tilde{\psi}(y(\tau)) < 0, \label{eq:DifferentialTimeForezen_jump} 
%		\end{cases} 
%	\end{align}
%\end{subequations}
\begin{subequations}\label{eq:DifferentialTimeForezen}
\begin{empheq}[left={{y}'(\tau) =\empheqlbrace}]{align}
\tilde{f}(y(\tau)),  \ &\tilde{\psi}(y(\tau)) \geq 0,\label{eq:DifferentialTimeForezen_ode}\\
	\tilde{\varphi}(y(\tau)) ,  \ &\tilde{\psi}(y(\tau)) < 0, \label{eq:DifferentialTimeForezen_jump} 
\end{empheq}
\end{subequations}
\\
with $\tilde{f}(y(\tau)) \coloneqq  (f(\tilde{x}(\tau)),1)$, ${\tilde{\varphi}(y(\tau)) \coloneqq (\varphi(\tilde{x}(\tau)),0)}$ and ${\tilde{\psi}(y(\tau)) \coloneqq \psi(\tilde{x}(\tau))}$. It is assumed that Assumption \ref{ass:Auxiliary_Dynamics} is satisfied.
\end{definition}
In the next section we illustrate the ideas and terminology on the example of the bouncing ball and provide some examples how to fulfill Assumption \ref{ass:Auxiliary_Dynamics}. 

\section{An Illustrating Example} \label{sec:TimeFreezingExample}
We consider the dynamics of a ball bouncing on a table given by \eqref{eq:bouncing_ball_model}. To mimic the restitution law, whenever ${\tilde{q}(\tau)< 0}$, we use the following linear ODE for the time interval $(\tau_{\textrm{s}},\tau_{\textrm{r}})$:
%\begin{subequations}
	\begin{align}\label{eq:restituion_dynamics}
		\tilde{q}'(\tau) &= \tilde{v}(\tau), \
	\tilde{v}'(\tau) = -k\tilde{q}(\tau)-c\tilde{v}(\tau), \
	{t}'(\tau) = 0, 
	\end{align}
%\end{subequations}
%here with $\tau_{\textrm{r}}$ we denote the \textit{restitution time} and the interval $(\tau_{\textrm{s}}^-,\tau_{\textrm{r}})$, when the ODE \eqref{eq:restituion_dynamics} is active as \textit{restitution phase}. 
The initial values are $\tilde{q}(\tau_{\textrm{s}}) = 0$, ${t}(\tau_{\textrm{s}}) = \tau_{\textrm{s}}$ and $\tilde{v}(\tau_{\textrm{s}})$ has the value of $v(\cdot)$ corresponding to the solution of \eqref{eq:bouncing_ball_model} at $\tau_{\textrm{s}}$, $k,\ c \in \R$ are parameters. The first two equations in \eqref{eq:restituion_dynamics} are a second-order linear ODE and can be solved analytically. Using so-called spring-damper systems to model mechanic impact is an old idea, cf. Chapter 2 in \cite{Brogliato2016}. \textcolor{black}{However, to recover the rigid-body impact dynamics as in \eqref{eq:bouncing_ball_model}, the system needs to get infinitely stiff, which makes it impractical in numerical computations \cite[Sec 2.4 ]{Brogliato2016}. Moreover, spring-damper models can cause negative contact forces \cite[Remark 2.3]{Brogliato2016}. Our approach does not suffer from these difficulties. The key difference here is the introduction of the clock state with time-freezing. As we will see below, this enables us to use even rather small values for $k$ to recover the exact impact law \eqref{eq:restitution_law}. We discard all pieces of the trajectory which correspond to the time-evolution of \eqref{eq:restituion_dynamics} and use just its end points, therefore the difficulties coming from standard compliant models are not part of the final trajectory of a time-frozen dynamic system.}
Since $\tilde{v}(\tau_{\textrm{s}})<0$ and $\tilde{q}(\tau_{\textrm{s}}) = 0$, with the right choice of the parameters $k$ and $c$ in \eqref{eq:restituion_dynamics}, we have ${\tilde{v}(\tau_{\textrm{r}})>0}$ and ${\tilde{q}(\tau_{\textrm{r}}) = 0}$. Afterwards we switch back to the dynamic system defined for $\tilde{q}(\tau)\geq 0$, which is discussed below. For the solution of \eqref{eq:restituion_dynamics} we require it to satisfy
% $v(\tau_{\textrm{r}}) = -\gamma v(\tau_{\textrm{s}}), y(\tau_{\textrm{r}}) = 0$.
\begin{align} \label{eq:restitution_conditions}
\tilde{v}(\tau_{\textrm{r}}) &= -\gamma \tilde{v}(\tau_{\textrm{s}}),\quad  \tilde{q}(\tau_{\textrm{r}}) = 0.
\end{align}
If $\gamma =1$ we simply pick some $k>0$ and set $c=0$. In the case $\gamma\in (0,1)$, using the analytic solution of the ODE \eqref{eq:restituion_dynamics} and assuming $c^2-4k<0$ we can select $k$ and $c$ so that the conditions \eqref{eq:restitution_conditions} are satisfied. For a fixed $k>0$ we can easily derive the following formula for $c$
%\begin{align}\label{eq:restitution_c}
%c &= 2|\ln(\gamma)|\sqrt{\frac{k}{\ln(\gamma)^2+\pi^2}}.
%\end{align}
\begin{align}\label{eq:restitution_c}
c &= 2|\ln(\gamma)|\sqrt{\nicefrac{k}{(\ln(\gamma)^2+\pi^2)}}.
\end{align}

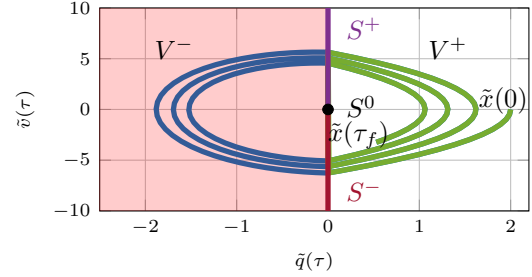
\begin{figure}[t]
	\centering
	\centering
	\vspace{0.1cm}
%	\tikzsetnextfilename{phase_plot}
	{\input{img/phase_plot4.tikz}}
	\caption{The state space and phase plot of the time-frozen dynamic system \eqref{eq:convex_comb_di}. The red shaded area $V^-$ is the prohibited region where the auxiliary dynamic flows and mimics the state jump so that the start point of every arc is in $S^-$ and end point in $S^+$. \textcolor{black}{The green curve shows the resulting trajectory after discarding the "time-frozen parts", cf. Theorem \ref{thm:recover_solution}.}}
	\vspace{-0.55cm}
	\label{fig:phase_plot}
\end{figure}
As already mentioned, after $\tau_{\textrm{r}}$ we switch back to the initial model with some modifications. Since $q(\tau_{\textrm{r}})>0$ (ball not in contact anymore) we can locally ignore equation \eqref{eq:restitution_law}. Furthermore, we add the dynamics of the clock state ${t'(\tau)=1}$, thus we get the following ODE:
%\begin{subequations}
	\begin{align}\label{eq:free_flight}
	\vspace{-0.5cm}
	\tilde{q}'(\tau) &= \tilde{v}(t), \
	m\tilde{v}'(\tau)= -mg, \
	{t}'(\tau) = 1,
	\end{align}
%\end{subequations}
%Note that using appropriate initial values, the same equations can be used for the pre-impact phase. 
Using the compact notation ${{y}(\tau)\coloneqq (\tilde{q}(\tau),\tilde{v}(\tau),t(\tau))}$ and denoting the r.h.s. of \eqref{eq:free_flight} in compact form as $f_1({y}(\tau))$ and analogously the r.h.s. of \eqref{eq:restituion_dynamics} as $f_2({y}(\tau))$, and defining $\tilde{\psi}({y}(\tau))\coloneqq\tilde{q}(\tau)$, we can write the combined dynamics in compact form as
\begin{align}\label{eq:convex_comb_di}
 {y}'(\tau) \in f_1({y})\alpha(\tilde{\psi}({y})) + f_2({y})(1-\alpha(\tilde{\psi}({y}))),
 \end{align}
where $\alpha(z)$ is a set-valued \textcolor{black}{step} function such that $\alpha(z) = 1$ if $z>0$, $\alpha(z) = 0$ if $z<0$, and $\alpha(z) \in [0,1]$ if $z=0$. The inclusion in \eqref{eq:convex_comb_di} accounts for the case if the dynamics stays on the manifold $\tilde{\psi}(y(\tau)) = 0$, which does not happen in our case. For more details see the concept of Filippov inclusions \cite{Filippov2013}. 
The equation \eqref{eq:convex_comb_di}  is an example of Definition \ref{def:ODE_time_frozeon} and the auxiliary dynamic system $y'(\tau) = f_2(y(\tau))$ satisfies all conditions of  Assumption \ref{ass:Auxiliary_Dynamics} by construction. Figure \ref{fig:phase_plot} depicts the phase plot of \eqref{eq:convex_comb_di}.
% note that in the \textit{prohibited region} $V^-$ the flow of the auxiliary dynamics mimics the state restitution law.

The set-valued \textcolor{black}{step} function $\alpha(z)$ can be represented as the solution of a parametric linear program (LP) \cite{Baumrucker2009}
\begin{align}\label{eq:indicator_lp}
\alpha(z) &=\underset{w}{\mathrm{argmin}}
-zw \textrm{ s.t. } 0\leq w\leq 1.
\end{align}
% or
%\vspace{-0.2cm}
Using the KKT conditions of this LP combined with \eqref{eq:convex_comb_di}, we obtain a DCS.
% We use the latter form for computational reasons, which is discussed in more detail in Section \ref{sec:NumericalExamples}.

%\sout{Note that for this formulation we can apply a number of different numerical simulation methods \cite{Acary2008}, which are not suitable for dynamics with state jumps.}
Observe that we got rid of the conditional algebraic restitution law \eqref{eq:restitution_law}. We have in fact a simpler nonsmooth dynamic system than \eqref{eq:bouncing_ball_model}, since the solution of \eqref{eq:convex_comb_di} is AC \cite{Filippov2013} and contains no jumps. We have in fact reduced the difficult case 3 with jumps to the simpler case 2 without jumps. Hence, there is no need to use measures, which simplifies the theoretical analysis as well as the numerical computation.
\begin{figure}[t]
	\centering
	\centering
	\vspace{0.3cm}
%	\tikzsetnextfilename{time_explain}
	{\input{img/time_explain2.tikz}}
	\caption{The clock state $t(\tau)$, with an illustration of the pseudo time, virtual time and physical time. The length of the pseudo time intervals is always the same $\tau_{\textrm{jump}}$.}
%	\vspace{-0.2cm}
	\label{fig:time_explain}
	%\end{figure}
	%
	%\begin{figure} [t]
	\centering
	\centering
	\vspace{-0.1cm}
	\tikzsetnextfilename{projected}
	{\input{img/projected2.tikz}}
	\caption{The velocity $v(\tau)$ and position $q(\tau)$ of the bouncing ball in pseudo time $\tau$ (top), and physical time $t(\tau)$ (bottom). The red shaded area in the top plot marks the intervals where the auxiliary dynamic system is active.}
	\vspace{-0.5cm}
	\label{fig:projected_solution}
\end{figure}

After getting rid of the state jump, the question is how to recover the true solution with state jumps? For illustration, we simulate \eqref{eq:convex_comb_di} with $\tilde{q}(0) = 10$, $\tilde{v}(0) = 0$ and $t(0) = 0$. We take $\gamma = 0.9$, where for a fixed $k= 20$, we obtain $c = 0.2998$ via \eqref{eq:restitution_c}. Figure  \ref{fig:time_explain} depicts the evolution of the clock state $t(\tau)$. We distinguish between three time concepts: 1) the \textit{pseudo time} $\tau$, which is the time of the nonsmooth dynamic system; 2) the \textit{physical time} $t$, the part of the pseudo time whenever ${t}'(\tau)>0$; 3) the \textit{virtual time} $t_\textrm{V}$, the part of the pseudo time $\tau$ whenever ${t}'(\tau)=0$ (restitution phases).
The top plot in Figure \ref{fig:projected_solution} depicts the state trajectories $\tilde{q}(\tau)$ and $\tilde{v}(\tau)$ in pseudo time $\tau$ and the bottom plot show the state trajectories in physical time $q(t(\tau))$ and $v(t(\tau))$. Obviously, we recover the true trajectories of the model in \eqref{eq:bouncing_ball_model}. The formal proof for this observation in a more general setting is provided in Section \ref{sec:Equivalence}.
%\vspace{-0.05cm}
%\sout{Note that if we pick a large $k$, the system gets more stiff, hence we would need a small step size for accurate simulations, on the other hand this makes the transition faster, so we can select a smaller simulation time $\tau_f$ to get to the desired final physical time $t(\tau_f)$. In fact, }
Using the analytic solution of \eqref{eq:restituion_dynamics} and equation \eqref{eq:restitution_c}, the length of the \textit{restitution phase} $\tau_{\textrm{jump}}$ can found to be
\vspace{-0.55cm}
\begin{align}\label{eq:T_rest}
 \tau_{\textrm{jump}} &=\sqrt{\nicefrac{(\pi^2+\ln(\gamma)^2)}{k}}.
\end{align}
%\vspace{-0.15cm}
% \sout{Moreover, it is easy to see that in this example with our model we do not obtain Zeno behavior (infinite switches in finite time), since every switch requires some virtual time to perform the state jump, which is upper bounded by the total pseudo simulation time $\tau_f$. }
\textcolor{black}{The dynamic system from Definition \ref{def:ODE_state_jump} excludes motion on the manifold $S$  \cite{Brogliato2016}, i.e. inelastic and persistent contacts ($\gamma =0$). The analysis of this case is different and beyond the scope of this paper.}

\textcolor{black}{
The difficult part in constructing a system from Definition \ref{def:ODE_time_frozeon} is to fulfill Assumption \ref{ass:Auxiliary_Dynamics}. In general any kind of compliant model can be used as long as its initial and final point satisfy the conditions in Assumption \ref{ass:Auxiliary_Dynamics}. We discuss briefly how to construct such systems for mechanical impact problems. Let $q\in \R^{m}$ be the generalized coordinates and $v \in R^{m}$ the generalized velocity of a rigid body.
%with $n_i\in\R^m,\ b\in\R$
Consider an affine unilateral constraint $\psi_i(q) = n_i^\top q + b_i \geq 0$. If the body collides with this constraint, then according to Newton's restitution law the post-impact velocity is: $v(t^+) = -\gamma n_i^\top v(t^-).$ As the velocity change happens only along the normal, we can project the system on the normal and perform the state jump law with the spring-damper model along this line and add the result back to the normal. This provides an auxiliary dynamic system satisfying Assumption \ref{ass:Auxiliary_Dynamics}:
\begin{align}\label{eq:impact_dynamics}
{\varphi}_i(\tilde{x}) &= N_i K N_i^\top(\tilde{q},\tilde{v}), \ N_i  \coloneqq \begin{bmatrix}
n_i & 0_{m,1} \\  0_{m,1} & n_i 
\end{bmatrix},
\end{align}
where $K \coloneqq \begin{bmatrix}
0 & 1 \\  -k & -c
\end{bmatrix}$ defines the two-dimensional linear spring-damper dynamics. Therefore, we can use for every affine constraint $\psi_i(q)$ an auxiliary dynamic of this form. For multiple constraints the Filippov convexification via step functions $\alpha(\cdot)$ (as generalization of \eqref{eq:convex_comb_di}) can be written using equation (4.1) in \cite{Dieci2011}. In case of activation of multiple perfect fiction-less constraints the negative reaction force is in the normal cone to the feasible set at this point \cite{Acary2008}, and hence the auxiliary dynamic has to evolve in this part of the state space. In case the constraints are orthogonal, the desired vector field is simply the sum of the neighboring fields, otherwise the analysis is a bit more involved.}

\section{Solution Relationship} \label{sec:Equivalence}
In this section we show how the solutions of the initial nonsmooth differential equation with a state jump law \eqref{eq:DifferentialJump} and the corresponding time-frozen system \eqref{eq:DifferentialTimeForezen} are related. Note that the function $t(\tau;x_0)$ is monotone by construction, e.g. Figure \ref{fig:time_explain}. Using the definitions from Section \ref{sec:TimeFreezing} we can state the main theoretical result.
\begin{theorem} \label{thm:recover_solution}
	Suppose that Assumption \ref{ass:Auxiliary_Dynamics} holds.
	Consider the initial value problem (IVP) corresponding to Definition \ref{def:ODE_time_frozeon} with a given $y_0 = (x_0,0)$ and ${x_0\in V^+}$ on a time interval $(0,\tau_{\textrm{f}})$, and the IVP corresponding to Definition \ref{def:ODE_state_jump} with the initial value $x_0$ on a time interval $(0,t_{\textrm{f}})\coloneqq(0,t(\tau_{\textrm{f}}))$, with $\tilde{x}(\tau_{\textrm{f}}) \in V^+$. Suppose that we have at most one time point $t_s$ where $\psi(x(t_s)) = 0$ on the time interval $(0,t_{\textrm{f}})$. Then the solution of the two IVPs: $x(t;x_0)$ and $y(\tau;y_0)$ fulfill at any $t \neq t_s$
	\begin{align}\label{eq:SolutionRelation}
	x(t(\tau)) = R y(\tau), 
\text{ with }
	R = \begin{bmatrix}
	\mathds{1}_{n_x,n_x}  &{0}_{n_x,1} \\ 
	{0}_{1,n_x} & 0 
	\end{bmatrix}.
	\end{align}
\end{theorem} 
\emph{Proof.}:
Denote the solution of IVP given by \eqref{eq:DifferentialTimeForezen_ode} and $y_0$ as $y_1(\tau;y_0)$ for some $\tau \in (0,\hat{\tau})$. Similarly, for \eqref{eq:DifferentialJump_ode} and ${x_0 \in V^+}$ for some $t(\tau) \in (0,t(\hat{\tau}))$ as  $x_1(t(\tau);x_0)$.  Note that if there is no $t_s \in (0,t_{\textrm{f}})$ such that $\psi(x(t_s)) = 0$ on this interval, then $t(\tau) = \int_{0}^{\tau} \dd \tau_1 = \tau$. Then setting $\hat{\tau} = \tau_{\textrm{f}}$, it follows that \eqref{eq:SolutionRelation} holds, since due to Definitions \ref{def:ODE_state_jump} and \ref{def:ODE_time_frozeon} it follows $x_1(t;x_0)= x(t;x_0)$ and  $y_1(\tau;x_0)= y(\tau;y_0)$.

If we have some  $t_s \in (0,t_{\textrm{f}})$ so that $\psi(x(t_s)) = 0$, then from the first part of the proof we have that \eqref{eq:SolutionRelation} holds for all $\tau \in(0,\tau_{\textrm{s}}^-)$ and hence for all $t(\tau) \in(0,t_s^-)$, with $t_s = \tau_{\textrm{s}}$. 
Its only left to prove that  \eqref{eq:SolutionRelation} holds for $\tau \in (\tau_{\textrm{s}}^+,\tau_{\textrm{f}})$ and the respective $t(\tau)$. Due to Assumption \ref{ass:Auxiliary_Dynamics} there exists a dynamic system $\tilde{x}'(\tau) = \varphi(\tilde{x}(\tau))$ that satisfies the restitution law and we have that $y(\tau_{\textrm{r}};y_0) = \Gamma(y(\tau_{\textrm{s}};y_0)) = y_s$. Note that $t'(\tau) = 0$ with $\tau \in (\tau_{\textrm{s}},\tau_{\textrm{r}})$, hence $t(\tau_{\textrm{r}}) = t(\tau_{\textrm{s}}) = t_s$.
Using this we have $y_1(\tau-\tau_{\textrm{r}},y_s) = y(\tau,y_0)$ for $\tau \in(\tau_{\textrm{r}},\tau_{\textrm{f}})$ and with denoting $x_s = Ry_s$, we see that $x_1(t(\tau)-t_s;x_s) = x(t(\tau),x_0)$ for $t(\tau) \in (t_s^+,t_{\textrm{f}})$.  
 Since the intervals $(t_s,t_{\textrm{f}})$ and $(\tau_{\textrm{r}},\tau_{\textrm{f}})$ have the same length and $x_s = Ry_s$, from the definitions of the corresponding IVPs, we conclude that relation \eqref{eq:SolutionRelation} holds. This completes the proof. \qed

%Obviously, the time-concatenation of $x_1(t_1,x_0)$ and $x_2(t_1;x_s)$ for $t_1\in (0,t_s^-)$ and $t_2 \in (t_s^+,t_f)$ is the solution $x(t,x_0)$, $t \in (0,t_f)$ of the IVP in Definition  \ref{def:ODE_state_jump}. 

The assumptions that we have at most one state jump on the time interval $(0,t_{\textrm{f}})$ can be always satisfied by shortening the regarded time interval and simplifies the proof without loss of generality. Furthermore, we avoid the analysis of the case with infinite switches in finite time (Zeno behavior).
For a desired physical simulation $t_{\textrm{f}}$ we always have to take a longer pseudo simulation time $\tau_{\textrm{f}} = t_{\textrm{f}} + N_{\textrm{J}}\tau_{\textrm{jump}}$ where $N_{\textrm{J}}$ is the number of state jumps on $(0,t_{\textrm{f}})$ for the original system. Obviously, we do not know a priori the number $N_{\textrm{J}}$. However, in OCPs this can be easily overcome with the use of a time transformation, which is shown in the next section. 

\section{Numerical Examples}\label{sec:NumericalExamples}
\subsection{Numerical Simulation}
\begin{figure}[t]
	\centering
	\centering
	\vspace{0.15cm}
	\tikzsetnextfilename{integrator_errors}
	{\input{img/numerical_k5.tikz}}
	\caption{\textcolor{black}{Accuracy of the explicit Euler scheme and Runge-Kutta 4 for different number of function evaluations for the time-freezing reformulation.}}
	\vspace{-0.55cm}
	\label{fig:numerical_error}
\end{figure}
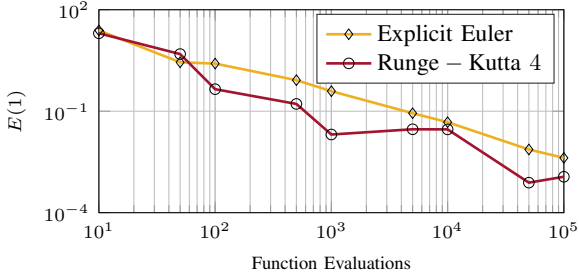

We first demonstrate the ease of use of the time-freezing in simulation problems. Consider again the example from Section \ref{sec:TimeFreezingExample}. The initial value is set to $x(0) = (0.5,0)$ and we simulate the original system for $t_{\textrm{f}}= 1$ s. For the time-freezing reformulation we take the compact form of \eqref{eq:convex_comb_di} and use the \textcolor{black}{step} function $ \alpha(z) = (1+\textrm{sign}(z))/2$. In this example the restitution coefficient is picked to be $\gamma = 0.9$, and for \textcolor{black}{$k =5$} we calculate \textcolor{black}{$c =0.1499$} via \eqref{eq:restitution_c}. The analytical solution has two jumps during the considered time interval. The pseudo simulation time is set to ${\tau_{\textrm{f}} = 1 + 2\tau_{\textrm{jump}}}$, where
$\tau_{\textrm{jump}} = \textcolor{black}{1.4058}$ and is obtained with equation \eqref{eq:T_rest}. In this numerical experiment we use the explicit Euler and Runge-Kutta 4 (RK4) schemes with equidistant steps. %Since we simulate a nonsmooth differential equation, the accuracy of both methods is $O(h)$ \cite{Acary2008}. 
The terminal numerical error, denoted as $E(1) = \| x(1) - \tilde{x}(t(\tau_{\textrm{f}}))\|_2$ is plotted over the number of function evaluations $M$ in the integrator, see Figure \ref{fig:numerical_error}.
%Thereby for some step size $h$, to advance the integrator the explicit Euler need one function evaluation and RK4 four. 
We clearly see that the error decreases for both methods with a smaller step size and that the numerical time-freezing solution converges to the analytic solution, which also confirms the result of Theorem 1. Due to the remaining nonsmoothness (case 2), the RK4 method does only achieve an order of one, as the Euler method. Opposed to standard spring-damper impact models, the most notable observation for our reformulation is that we do not need a large $k$ (which makes the system very stiff and costly to integrate) to get a very accurate numerical approximation.
%\begin{figure}[t]
%	\centering
%	\centering
%	\vspace{0.2cm}
%	\tikzsetnextfilename{integrator_errors}
%	{\input{img/integrator_errors.tikz}}
%	\caption{Accuracy of the explicit Euler scheme and Runge-Kutta 4 for different number of function evaluations for the time-freezing reformulation.}
%	\vspace{-0.5cm}
%	\label{fig:numerical_error}
%\end{figure}
\subsection{Numerical Optimal Control}
\color{black}
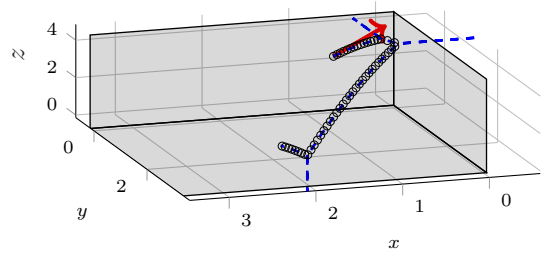
\begin{figure}[t]
	\centering
	\centering
	\vspace{0.1cm}
	%	\tikzsetnextfilename{3d_simulation}
	{\input{img/sim3d.tikz}}
	\caption{\textcolor{black}{Simulation of the particle trajectory with the given initial conditions and  $u=0_{3,1}$. The black dashed curve is the position $q(\tau)$ in pseudo time, the black circles show the particle in physical time $q(t)$.}}
	\vspace{-0.6cm}
	\label{fig:3d_simulation}
\end{figure}
%\color{black}
\color{black}
We consider a time-optimal control problem  of a moving 3D particle. The particle is represented via its position $q\coloneqq(q_x,q_y,q_z)$ and velocity $v \coloneqq (v_x,v_y,v_z)$. The mass of the particle is $m= 1$ kg and it is controlled via a bounded magnetic force $u \coloneqq (F_x,F_y,F_z)$. The particle's initial position is $q(0) = (4,4,1)$ and the initial velocity is $v(0) = (-3,-3.5,0)$.
The free flight dynamics of the particle with the clock state are given by $ {y}' = \tilde{f}(y,u) \coloneqq (\tilde{v}_x,\tilde{v}_y,\tilde{v}_y,F_x/m,F_x/m,(F_y-g)/m,1)$. We have three unilateral constraints given by 
$\psi_1(x) = (1,0,0)^\top q$, $\psi_2(y) =  (0,1,0)^\top q$ and $\psi_3(y) =  (0,0,1)^\top q$. For every constraint we define an auxiliary dynamic system according to equation \eqref{eq:impact_dynamics}. Since the angle between every two constraints is $\pi/2$, the auxiliary dynamics for multiple active constraint at the corners (evolving in the normal cones of the feasible set at corners) is the sum of the vector fields of each active constraint. We exploit this to simplify the Filippov representation (see equation (4.1.) in \cite{Dieci2011}) via step functions and avoid all possible combinations. This yields the following dynamics
\begin{align}
\begin{split}
y'&\in \prod_{i=1}^3 \alpha(\tilde{\psi_i}(y)) \tilde{f}(y,u) + \sum_{i =1}^{3}(1-\alpha (\tilde{\psi_i}(y)) \tilde{\varphi}_i(y) 
\end{split}
\end{align}
The r.h.s. of the last equation is compactly denoted as $F(y,u)$. The trajectory of the unactuated particle ($u=0_{3,1}$) is depicted in Figure \ref{fig:3d_simulation}. The particle firsts hits the wall in the $x$-$z$ plane, then the wall in the $y$-$z$ plane and moves away from the corner. The goal in the OCP is to have the particle in minimum time at the final position $q_\textrm{traget} = (5,5,1)$, on the same line connecting the corner and $q(0)$ in the $x$-$y$ plane. The time-optimal control problem reads as
\begin{mini!}[2]
	{y(\cdot),u(\cdot),w}
	{  t(1)  + \rho || q(1) -q_{\textrm{target}} ||_2^2 }
	{\label{eq:OCP_time_freezing}}{}
	\addConstraint{{y}(0) - y_0 = 0}{\label{eq:OCP_time_freezing_IV}}
	\addConstraint{{y}'(s)\in w F(y(s),u(s)), \ s\in [0,1]}{\label{eq:OCP_time_freezing_ODE1}}
	\addConstraint{-mg e_3 \leq u(s) \leq mg e_3, \ s\in [0,1]}{\label{eq:OCP_time_freezing_ForceLimit}}
	\addConstraint{w_{\textrm{max}}^{-1} \leq w(s) \leq  w_{\textrm{max}} , \ s\in [0,1].}{\label{eq:OCP_time_freezing_SpeedOfTime}}
\end{mini!}
\noindent \hspace{-0.38cm}where $y_0 = (q(0),v(0),0)$. To achieve a time-optimal formulation, in this OCP another time transformation $\tau = ws$ is used, where $w$ is a parameter (the "speed of pseudo-time") and $s$ is the new pseudo-time. Moreover, we also have an upper and lower bound on $w\in [w_{\textrm{max}}^{-1},w_{\textrm{max}}]$ with $w_{\textrm{max}}=20$, to avoid numerical difficulties. 
All time-derivatives are now w.r.t. $s$, hence all differential equations are scaled by $w$. For an initial guess for $w$ we set $2T$ (with $T=1$) as it is likely to be greater than $T$ since the auxiliary dynamics take some of the "time budget". The penalty parameter $\rho$ is set to $10^2$. In the dynamics \eqref{eq:OCP_time_freezing_ODE1} we replace all step functions $\alpha(\tilde{\psi}_i(y(s))),\ i = 1,2,3$, with the KKT conditions of the parametric LP formulation \eqref{eq:indicator_lp} and obtain a DCP.
\color{black}
\begin{figure}[t]
	\centering
	\centering
	\vspace{0.1cm}
	%	\tikzsetnextfilename{3d_simulation}
	{\input{img/ocp3d.tikz}}
	\caption{\textcolor{black}{Resulting trajectory after solving the time-optimal control problem with final time $t(1)= 0.87$ s. The black dashed curve is the position $q(\tau)$ in pseudo time, the black circles show the particle in physical time $q(t)$.}}
	\vspace{-0.62cm}
	\label{fig:3d_ocp}
\end{figure}
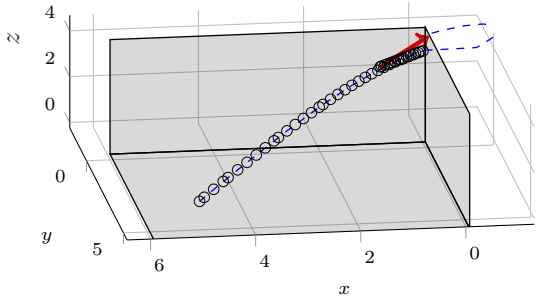
\color{black}
\color{black}
We use a fully simultaneous approach and discretize the infinite-dimensional OCP \eqref{eq:OCP_time_freezing} using the implicit Euler scheme with a step size ${h=0.005}$. The discretized control inputs are taken to be constant over the finite elements. Since we discretize a DCP, the discretized OCP yields an NLP which is an MPCC. To solve this MPCC we use a  homotopy penalization approach, which works as follows. The complementarity constraints are penalized with a positive parameter $\mu >0$ and added to the objective  (e.g.  $a^\top b = 0$ is added as $\mu a^\top b$), hence we obtain a smooth NLP.
Furthermore, if $\mu$ is  larger than a critical value of the penalty parameter, then the CCs will be satisfied at the solution \cite{Ralph2004}, which means one needs only to solve a single NLP.
To solve an MPCC originating from an OCP, special care has to be taken. In \cite{Nurkanovic2020} it was shown that in discretized OCPs resulting in MPCCs, we have to take a step size $h=o(\sigma)$, where $\sigma$ is a relaxation parameter for the complementarity conditions, e.g. $a^\top b \leq \sigma$. This is needed to obtain the right numerical sensitivities and avoid getting stuck in spurious local solutions close to the initial guess. In general, a one-to-one relation between KKT points of relaxation and penalization schemes can be established \cite{Leyffer2006}.   Therefore, we solve a sequence of NLPs for a varying penalty parameter to avoid convergence to spurious solutions. The parameter $\mu$ is updated by the following rule: $\mu_{i+1} = 10 \mu_{i}$ with $\mu_0 = 10^{-3}$, where $i$ is the number of the problem in the sequence. The primal solution of the previous problem is used as a solution guess for the next problem in the sequence and we solve in this example in total 7 problems with \texttt{IPOPT} \cite{Waechter2006} via its \texttt{CasADi} \cite{Andersson2018} interface. 
The solution trajectory $q(\cdot)$ of the OCP is given in Figure \ref{fig:3d_ocp}. The algorithm finds a trajectory with multiple simultaneous impacts. This is the solution one would intuitively expect, since the target point lies on the line connecting the corner and initial point. There is no need for integer variables, nor for a good solution guess, nor the need to incorporate the algebraic restitution law explicitly (it is not even clear how this could be done in a smooth optimization problem formulation). Moreover, the CC are satisfied at the solution, hence no "smoothing effects" are left at the solution. 
%\sout{Observe that if many switches occur, $w$ has to take a larger value a terminal physical time $T$ can be reached, since every restitution phase "uses" a part of the available pseudo time. This might require one to use a smaller step size $h$ to achieve the desired accuracy.}
\color{black}
%\vspace{-0.2cm}
\section{Conclusions and Outlook}\label{sec:Conclusions}
In this paper we proposed a novel reformulation for differential equations with state jumps into a significantly easier class of problems. We also provide a proof that the solutions of the two systems are related and how to recover the solution of the original system. The proposed reformulation significantly simplifies to use differential equations with state jumps in numerical optimal control. 
The efficacy of the approach is illustrated on a simulation example and \textcolor{black}{a time-optimal control problem where we obtain a solution with multiple simultaneous impacts. The hard part, in general, is to find an auxiliary dynamic system satisfying Assumption \ref{ass:Auxiliary_Dynamics}. Depending on the application, a good starting point are existing compliant models. A systematic way to obtain such  differential equations is subject of future research.}
%\color{black}
%\begin{figure}[t]
%	\centering
%	\centering
%	\vspace{0.2cm}
%	%	\tikzsetnextfilename{3d_simulation}
%	{\input{img/ocp3d.tikz}}
%	\caption{\textcolor{blue}{Resulting trajectory after solving the time-optimal control problem with final time $t(1)= 0.87$ s. The blue dashed curve is the position $q(\tau)$ in pseudo time, the black circles show the particle in physical time $q(t)$.}}
%	\vspace{-0.6cm}
%	\label{fig:3d_ocp}
%\end{figure}
%\color{black}
%An interesting practical case where this approach might be extend is for rigid-body with impacts and inelastic collisions where the contact is not immediately broken, since they form a rich source of OCPs.

\bibliographystyle{ieeetran}
\bibliography{bib/syscop}
\end{document}

%% file: img/phase_plot4.tikz
\definecolor{mycolor1}{rgb}{0.00000,0.44700,0.74100}%
\definecolor{mycolor2}{rgb}{0.46600,0.67400,0.18800}%
\definecolor{mycolor3}{rgb}{0.49400,0.18400,0.55600}%
\definecolor{mycolor4}{rgb}{0.63500,0.07800,0.18400}%
\pgfplotsset{compat=1.13}
\setlength{\fwidth}{6.0cm}
\setlength{\fheight}{2.7cm}
%\pgfplotsset{every tick label/.append style={font=\tiny}} 
\pgfplotsset{
	every axis plot/.append style={line width=1.0pt},
	every axis plot post/.append style={
		every mark/.append style={line width=0.0pt}
		%every tick label/.append style={font=\tiny}
	}
}

\begin{tikzpicture}

\begin{axis}[%
width=0.951\fwidth,
height=\fheight,
at={(0\fwidth,0\fheight)},
scale only axis,
xmin=-2.5,
xmax=2.2,
xlabel style={font=\color{white!15!black}},
xlabel={$\tilde{q}(\tau)$},
ymin=-10,
ymax=10,
ylabel style={font=\color{white!15!black}},
ylabel={$\tilde{v}(\tau)$},
axis background/.style={fill=white},
xmajorgrids,
ymajorgrids,
legend style={legend cell align=left, align=left, draw=white!15!black},
xlabel style={font={\scriptsize}},ylabel style={font=\scriptsize},ylabel shift={-0cm},ticklabel style={font=\scriptsize}
]
\addplot [color=mycolor1, line width=2.0pt]
  table[row sep=crcr]{%
2	0\\
1.9877375	-0.4905\\
1.95095	-0.981\\
1.8896375	-1.4715\\
1.8038	-1.962\\
1.690980095	-2.46231\\
1.55264438	-2.96262\\
1.388792855	-3.46293\\
1.195457375	-3.97305\\
0.975595655	-4.48317\\
0.724209500000001	-5.0031\\
0.445267055	-5.52303\\
0.132720455	-6.05277\\
-0.00907627344967121	-6.26219406503727\\
-0.1339536026088	-6.22140482372815\\
-0.263922109540663	-6.15205045267392\\
-0.398205760252651	-6.05070881314124\\
-0.535859894956697	-5.91404981084725\\
-0.675760834102242	-5.73891374850103\\
-0.816598226163694	-5.5224000669777\\
-0.962127519513847	-5.25125669857591\\
-1.10477407799674	-4.93144195749749\\
-1.24251461936882	-4.56162635927286\\
-1.37314934433735	-4.14126828393253\\
-1.49432461771154	-3.67074508612524\\
-1.60670721546618	-3.13476443275977\\
-1.70344595506184	-2.55096831550696\\
-1.77988779848492	-1.94200998020611\\
-1.83660462974947	-1.29644129019398\\
-1.87114034307719	-0.620996354741527\\
-1.88128189476803	0.0573222905826736\\
-1.86708395017449	0.729744702248231\\
-1.8289146590908	1.38764365297624\\
-1.76946254251912	2.00538719031203\\
-1.68887624147199	2.59411415585519\\
-1.59143375620293	3.13162799312168\\
-1.48000013495346	3.61523079636791\\
-1.35734687763629	4.04357071704452\\
-1.22610460112344	4.41648242399416\\
-1.08872690315628	4.73482333765136\\
-0.947464493755747	5.00031105861902\\
-0.809560404688878	5.20840165875763\\
-0.67193681537649	5.37191369124219\\
-0.541530747755249	5.48994176731969\\
-0.414251047699607	5.57291131901392\\
-0.291030571981911	5.62439164110509\\
-0.166978687841544	5.648484706848\\
-0.0426980353859641	5.64520165454433\\
0.00243334401027973	5.63352978073956\\
0.292340972170217	5.10378978073956\\
0.549063685549413	4.58385978073956\\
0.778230108928609	4.06392978073956\\
0.976291337527067	3.55380978073956\\
1.14782632612552	3.04368978073956\\
1.29029659994324	2.54337978073956\\
1.40725106376096	2.04306978073956\\
1.49868971757868	1.54275978073956\\
1.56356520661565	1.05225978073956\\
1.60391569565263	0.561759780739563\\
1.61974118468961	0.0712597807395596\\
1.61104167372659	-0.419240219260438\\
1.57781716276357	-0.909740219260439\\
1.52006765180055	-1.40024021926044\\
1.43779314083752	-1.89074021926044\\
1.32860748465524	-2.39105021926043\\
1.19390601847296	-2.89136021926044\\
1.03368874229068	-3.39167021926043\\
0.844058770889134	-3.90179021926044\\
0.627902559487592	-4.41191021926044\\
0.385220108086049	-4.92203021926044\\
0.110574361465246	-5.44196021926044\\
-0.00586394333942586	-5.63649192111383\\
-0.129464119776606	-5.59538734003365\\
-0.251864300898813	-5.52746628893488\\
-0.377910005160065	-5.42827817463208\\
-0.506647100447572	-5.29477634883703\\
-0.642069612646077	-5.1166346883014\\
-0.777322362151168	-4.89621071207588\\
-0.915418913422246	-4.62129553838276\\
-1.04930307275505	-4.29796269411656\\
-1.17686363157668	-3.92556710537428\\
-1.29933379732944	-3.49063739232967\\
-1.40987615401997	-3.00617110896486\\
-1.50588349711969	-2.47501830897239\\
-1.58473058122845	-1.90137336413492\\
-1.64383958958106	-1.29086075142446\\
-1.68009620778972	-0.667535340334605\\
-1.69321664987761	-0.0230487854027945\\
-1.68190928032125	0.616597992208384\\
-1.64653251608952	1.24225549133898\\
-1.5896249496198	1.82955098349986\\
-1.51150987378047	2.38733784882735\\
-1.41631846383669	2.89471547159462\\
-1.30693960442985	3.34856187817102\\
-1.18618751216791	3.74717940868946\\
-1.05674773373704	4.09013516217687\\
-0.925509062242856	4.36978776060733\\
-0.79087447082298	4.59887603837041\\
-0.659562962870988	4.7742152137777\\
-0.533690264485042	4.90272687125138\\
-0.409909633451695	4.99448448490329\\
-0.2892831899632	5.0528441242567\\
-0.172689245465792	5.08126681713874\\
-0.0557514996688075	5.08274860595113\\
0.0051848935788037	5.06395461946733\\
0.259796343410572	4.54402461946733\\
0.48685150324234	4.02409461946733\\
0.682841303454641	3.51397461946733\\
0.852304863666942	3.00385461946733\\
0.992743544259778	2.50354461946732\\
1.10766641485261	2.00323461946733\\
1.19707347544544	1.50292461946733\\
1.25995720641881	1.01242461946733\\
1.29831593739218	0.521924619467323\\
1.31214966836554	0.0314246194673338\\
1.30145839933891	-0.459075380532673\\
1.26624213031228	-0.949575380532672\\
1.20650086128564	-1.44007538053267\\
1.12223459225901	-1.93057538053267\\
1.01101734285184	-2.43088538053267\\
0.874284283444675	-2.93119538053267\\
0.712035414037508	-3.43150538053267\\
0.520334014249811	-3.94162538053267\\
0.30210637446211	-4.45174538053267\\
0.0573524946744124	-4.96186538053267\\
-0.00796724345320499	-5.07223742747132\\
-0.124221150805472	-5.03237360355826\\
-0.244218798907346	-4.96268914098382\\
-0.367067550659527	-4.86018233225124\\
-0.491705040556427	-4.72204225363031\\
-0.621435993900485	-4.53858105742615\\
-0.749861288778956	-4.31235959657773\\
-0.879300498430754	-4.0321753738036\\
-1.0068850944888	-3.69360086724475\\
-1.12597483649462	-3.30539488144276\\
-1.23412779920147	-2.86900407325924\\
-1.32882849737584	-2.38707580735505\\
-1.40939879919903	-1.8490850379484\\
-1.47035844622643	-1.27388032287858\\
-1.50857243989983	-0.68414015589163\\
-1.52373122969759	-0.0736517190574588\\
-1.51451774792048	0.532853805141809\\
-1.4812843840778	1.12575986332869\\
-1.42645153559704	1.68184940898752\\
-1.35049063974446	2.20791719120067\\
-1.25741505820178	2.68426158051148\\
-1.15013667347251	3.10739961646573\\
-1.03151067812389	3.4753461616168\\
-0.908058321060579	3.77919507357823\\
-0.783104123341113	4.02350061920687\\
-0.659448642136665	4.2138126492956\\
-0.535028249738764	4.36073899117344\\
-0.415821922293377	4.46389597378343\\
-0.29880813784601	4.53203059826469\\
-0.184993104851134	4.56839712762009\\
-0.0752033510464418	4.57635026764959\\
0.00252434813327618	4.56116401037674\\
0.230487895683243	4.04123401037674\\
0.427368944222833	3.53111401037675\\
0.597723752762425	3.02099401037675\\
0.73903654229164	2.52068401037674\\
0.854833521820853	2.02037401037674\\
0.945114691350066	1.52006401037675\\
1.0088553918689	1.02956401037674\\
1.04807109238774	0.539064010376743\\
1.06276179290658	0.0485640103767437\\
1.05292749342542	-0.441935989623254\\
1.01856819394425	-0.932435989623251\\
0.95968389446309	-1.42293598962326\\
0.876274594981928	-1.91343598962326\\
0.765931454511142	-2.41374598962326\\
0.630072504040357	-2.91405598962325\\
0.528567111382789	-3.23778598962326\\
};

\addplot [color=mycolor2, line width=2.0pt]
  table[row sep=crcr]{%
2	0\\
1.9877375	-0.4905\\
1.95095	-0.981\\
1.8896375	-1.4715\\
1.8038	-1.962\\
1.690980095	-2.46231\\
1.55264438	-2.96262\\
1.388792855	-3.46293\\
1.195457375	-3.97305\\
0.975595655	-4.48317\\
0.724209500000001	-5.0031\\
0.445267055	-5.52303\\
0.132720455	-6.05277\\
0.00970305500000013	-6.24897\\
nan	nan\\
0.00243334401027973	5.63352978073956\\
0.292340972170217	5.10378978073956\\
0.549063685549413	4.58385978073956\\
0.749542255463433	4.13259978073956\\
nan	nan\\
0.757787835024912	4.11297978073956\\
0.958399663623369	3.60285978073956\\
1.13248525222183	3.09273978073956\\
1.27745707603954	2.59242978073956\\
1.39691308985726	2.09211978073956\\
1.49085329367498	1.59180978073956\\
1.55818128271196	1.10130978073956\\
1.60098427174894	0.610809780739562\\
1.61926226078591	0.120309780739564\\
1.61301524982289	-0.370190219260438\\
1.58224323885987	-0.860690219260438\\
1.52694622789685	-1.35119021926044\\
1.44712421693383	-1.84169021926044\\
1.34044011075154	-2.34200021926044\\
1.20824019456926	-2.84231021926043\\
1.05052446838698	-3.34262021926044\\
0.867292932204697	-3.84293021926044\\
0.654197440803154	-4.35305021926044\\
0.414575709401613	-4.86317021926044\\
0.143049542780807	-5.38310021926044\\
0.00540641229929673	-5.62835021926044\\
0.000116033959336903	5.07376461946733\\
0.255247413791105	4.55383461946733\\
0.482822503622873	4.03390461946733\\
0.679322423835174	3.52378461946733\\
0.849296104047474	3.01366461946733\\
0.990235094640307	2.51335461946733\\
1.10565827523314	2.01304461946733\\
1.19556564582598	1.51273461946733\\
1.25893987679934	1.02223461946732\\
1.29778910777271	0.531734619467334\\
1.31211333874607	0.0412346194673283\\
1.30191256971944	-0.44926538053267\\
1.26718680069281	-0.939765380532666\\
1.20793603166617	-1.43026538053268\\
1.12416026263954	-1.92076538053266\\
1.01344332323238	-2.42107538053267\\
0.87721057382521	-2.92138538053266\\
0.715462014418043	-3.42169538053267\\
0.524270734630345	-3.93181538053267\\
0.306553214842644	-4.44193538053267\\
0.0623094550549421	-4.95205538053268\\
0.0021784704885528	-5.06977538053267\\
0.00252434813327618	4.56116401037674\\
0.230487895683243	4.04123401037674\\
0.427368944222833	3.53111401037675\\
0.597723752762425	3.02099401037675\\
0.73903654229164	2.52068401037674\\
0.854833521820853	2.02037401037674\\
0.945114691350066	1.52006401037675\\
1.0088553918689	1.02956401037674\\
1.04807109238774	0.539064010376743\\
1.06276179290658	0.0485640103767437\\
1.05292749342542	-0.441935989623254\\
1.01856819394425	-0.932435989623251\\
0.95968389446309	-1.42293598962326\\
0.876274594981928	-1.91343598962326\\
0.765931454511142	-2.41374598962326\\
0.630072504040357	-2.91405598962325\\
0.528567111382789	-3.23778598962326\\
};

\addplot[area legend, draw=none, fill=red, fill opacity=0.2]
table[row sep=crcr] {%
x	y\\
-2.5	-10\\
-2.5	10\\
0	10\\
0	-10\\
}--cycle;

\node[right, align=left]
at (axis cs:-2,6) {$V^-$};
\node[right, align=left]
at (axis cs:1,6) {$V^+$};
\addplot [color=mycolor3, line width=2.0pt]
  table[row sep=crcr]{%
0	0\\
0	10\\
};

\addplot [color=mycolor4, line width=2.0pt]
  table[row sep=crcr]{%
0	-0\\
0	-10\\
};

\addplot [color=black, line width=10.0pt, draw=none, mark=*, mark options={solid, black}]
  table[row sep=crcr]{%
0	0\\
};

\node[right, align=left, font=\color{mycolor3}]
at (axis cs:0.1,8) {$S^+$};
\node[right, align=left, font=\color{mycolor4}]
at (axis cs:0.1,-8) {$S^-$};
\node[right, align=left]
at (axis cs:0.1,0.3) {$S^0$};
\node[right, align=left]
at (axis cs:1.55,0.85) {$\tilde{x}(0)$};
\node[right, align=left]
at (axis cs:-0.1,-2.5) {$\tilde{x}(\tau_f)$};
\end{axis}

\begin{axis}[%
width=1.227\fwidth,
height=1.227\fheight,
at={(-0.16\fwidth,-0.135\fheight)},
scale only axis,
xmin=0,
xmax=1,
ymin=0,
ymax=1,
axis line style={draw=none},
ticks=none,
axis x line*=bottom,
axis y line*=left,
legend style={legend cell align=left, align=left, draw=white!15!black},
xlabel style={font={\scriptsize}},ylabel style={font=\scriptsize},ylabel shift={-0cm},ticklabel style={font=\scriptsize}
]
\end{axis}
\end{tikzpicture}%

%% file: img/time_explain2.tikz
\definecolor{mycolor1}{rgb}{0.00000,0.44700,0.74100}%
\definecolor{mycolor2}{rgb}{0.85000,0.32500,0.09800}%

\pgfplotsset{compat=1.13}
\setlength{\fwidth}{6.5cm}
\setlength{\fheight}{2.5cm}
%\pgfplotsset{every tick label/.append style={font=\tiny}} 
\pgfplotsset{
	every axis plot/.append style={line width=1.0pt},
	every axis plot post/.append style={
		every mark/.append style={line width=0.0pt}
		%every tick label/.append style={font=\tiny}
	}
}

\begin{tikzpicture}

\begin{axis}[%
width=0.951\fwidth,
height=\fheight,
at={(0\fwidth,0\fheight)},
scale only axis,
xmin=0,
xmax=15,
xmajorgrids,
ymajorgrids,
xlabel style={font=\color{white!15!black}},
xlabel={$\tau$ [pseudo time]},
ymin=0,
ymax=12,
ylabel style={font=\color{white!15!black}},
ylabel={$t(\tau)$ [physical time]},
axis background/.style={fill=white},
legend style={legend cell align=left, align=left, draw=white!15!black},
xlabel style={font={\scriptsize}},ylabel style={font=\scriptsize},ylabel shift={-0cm},ticklabel style={font=\scriptsize}
]
\addplot [color=mycolor1, line width=1.5pt]
  table[row sep=crcr]{%
0	0\\
1.495	1.42784312292855\\
2.18	1.47712357987742\\
4.77	3.99796074420031\\
5.45	4.04424715975103\\
7.785	6.31106660334263\\
8.47	6.36137073962435\\
10.57	8.39286187656759\\
11.25	8.438494319496\\
13.145	10.2664776224657\\
13.83	10.3156178993727\\
15	11.4856178993727\\
};
%\addlegendentry{$\text{data}1$}
      \node[anchor=west,scale=0.81] (source) at (axis cs:10,1.5){ Virtual Time};
      \node (destination) at (axis cs:1.6,1.3){};
      \draw[->](source)--(destination);
      \node (destination) at (axis cs:4.8,3.9){};
      \draw[->](source)--(destination);
      \node (destination) at (axis cs:7.8,6.5){};
       \draw[->](source)--(destination);
       \node (destination) at (axis cs:10.8,8.8){};
       \draw[->](source)--(destination);
       \node (destination) at (axis cs:13.5,10.6){};
       \draw[->](source)--(destination);
       
\end{axis}
\end{tikzpicture}%

%% file: img/projected2.tikz
\definecolor{mycolor1}{rgb}{0.00000,0.44700,0.74100}%
\definecolor{mycolor2}{rgb}{0.85000,0.32500,0.09800}%

\pgfplotsset{compat=1.13}
\setlength{\fwidth}{7.0cm}
\setlength{\fheight}{4.5cm}
%\pgfplotsset{every tick label/.append style={font=\tiny}} 
\pgfplotsset{
	every axis plot/.append style={line width=1.0pt},
	every axis plot post/.append style={
		every mark/.append style={line width=0.0pt}
		%every tick label/.append style={font=\tiny}
	}
}

\begin{tikzpicture}

%\centering
\begin{axis}[%
width=0.952\fwidth,
height=0.410\fheight,
at={(0\fwidth,0.582\fheight)},
scale only axis,
xmin=0,
xmax=15,
xlabel style={font=\color{white!15!black}},
xlabel={$\tau$ [pseudo time]},
ymin=-18,
ymax=15,
ylabel style={font=\color{white!15!black}},
ylabel={$\tilde{q}(\tau),\ \tilde{v}(\tau)$},
axis background/.style={fill=white},
xmajorgrids,
ymajorgrids,
legend style={at={(0.55,0.015)}, anchor=south west, legend cell align=left, align=left, draw=white!15!black, nodes={scale=0.70, transform shape}},
legend columns=2, 
xlabel style={font={\scriptsize}},ylabel style={font=\scriptsize},ylabel shift={-0cm},ticklabel style={font=\scriptsize}
]
\addplot [color=mycolor1]
  table[row sep=crcr]{%
0	10\\
0.130000000000001	9.9171055\\
0.279999999999999	9.615448\\
0.465	8.939416375\\
0.695	7.630762375\\
0.965	5.43234137500001\\
1.285	1.900741375\\
1.715	-2.87877107278431\\
1.78	-2.97269827785079\\
1.855	-2.77274053779275\\
1.96	-1.99983898893506\\
2.205	0.909347392292215\\
2.56	4.50779092292722\\
2.86	6.58489661783004\\
3.105	7.62625947700068\\
3.295	8.02844875043913\\
3.435	8.09818774139378\\
3.565	7.99077987585168\\
3.72	7.64603827655146\\
3.91	6.90192904998992\\
4.145	5.4916965526638\\
4.425	3.1041185345731\\
4.765	-0.790185486829671\\
4.96	-2.48550536333126\\
5.05	-2.67650484964091\\
5.125	-2.50812229401644\\
5.235	-1.78237184435916\\
5.495	0.99484158140756\\
5.825	3.90927436724586\\
6.105	5.54435242795715\\
6.325	6.28950661851601\\
6.495	6.54010608394787\\
6.625	6.54044596927811\\
6.76	6.36532255212105\\
6.925	5.9084853200402\\
7.13	4.96885541806097\\
7.38	3.26502154369605\\
7.675	0.46589619694544\\
7.97	-2.21197103657061\\
8.065	-2.40907828088944\\
8.145	-2.24193534637131\\
8.26	-1.53194100786328\\
9.025	4.38372697361484\\
9.245	5.08642286603081\\
9.41	5.30185466034278\\
9.54	5.28348205131585\\
9.68	5.07828716467146\\
9.85	4.5706284451747\\
10.06	3.55210161520813\\
10.315	1.73370866096299\\
10.825	-2.15949831235803\\
10.905	-2.09251123682439\\
11.005	-1.64768582184394\\
11.185	-0.17828412606562\\
11.51	2.35060631289597\\
11.77	3.62890436470737\\
11.97	4.16095055840844\\
12.12	4.30247270368426\\
12.25	4.24658322958996\\
12.395	3.98865809502324\\
12.575	3.38153266935421\\
12.795	2.2078504824254\\
13.06	0.16368181407932\\
13.325	-1.76946504339561\\
13.43	-1.95043186631955\\
13.515	-1.78657707860888\\
13.635	-1.14757413284609\\
14.25	2.8007410180141\\
14.455	3.34671456474522\\
14.605	3.48501615381677\\
14.735	3.42633553101212\\
14.88	3.16529719211462\\
15	2.79328646337186\\
};
\addlegendentry{$\tilde{q}(\tau)$}

\addplot [color=mycolor2]
  table[row sep=crcr]{%
0	0\\
1.43	-13.9974352153923\\
1.48	-13.4144368563016\\
1.57	-10.7629222755509\\
1.74	-1.89341488667716\\
2.01	11.2340969542043\\
2.115	12.6347411500197\\
2.16	12.3191856496761\\
4.7	-12.5982143503239\\
4.75	-12.1214810333755\\
4.835	-9.96591653684869\\
4.995	-2.66304167434867\\
5.285	10.1574624897695\\
5.39	11.3711372363155\\
5.45	10.8917145025403\\
7.72	-11.3338239242465\\
7.775	-10.7721260947426\\
7.87	-8.35415589913626\\
8.055	-0.27764402330774\\
8.3	9.12119668172544\\
8.405	10.2341836161548\\
8.455	9.86487223825436\\
10.505	-10.1992134706011\\
10.565	-9.61468947154309\\
10.665	-7.19637760986771\\
10.87	1.05963282169206\\
11.09	8.30727889642483\\
11.19	9.21073686568865\\
11.24	8.84053096850539\\
13.08	-9.18415695544899\\
13.135	-8.74055031025775\\
13.23	-6.79872162792052\\
13.415	-0.270389943267327\\
13.66	7.36828512572903\\
13.765	8.28963534106107\\
13.805	8.03426059381034\\
15	-3.68868940618966\\
};
\addlegendentry{$\tilde{v}(\tau)$}

\addplot[area legend, draw=none, fill=red, fill opacity=0.2, forget plot]
table[row sep=crcr] {%
x	y\\
1.425	-18\\
1.425	15\\
2.13	15\\
2.13	-18\\
}--cycle;

\addplot[area legend, draw=none, fill=red, fill opacity=0.2, forget plot]
table[row sep=crcr] {%
x	y\\
4.7	-18\\
4.7	15\\
5.4	15\\
5.4	-18\\
}--cycle;

\addplot[area legend, draw=none, fill=red, fill opacity=0.2, forget plot]
table[row sep=crcr] {%
x	y\\
7.715	-18\\
7.715	15\\
8.415	15\\
8.415	-18\\
}--cycle;

\addplot[area legend, draw=none, fill=red, fill opacity=0.2, forget plot]
table[row sep=crcr] {%
x	y\\
10.5	-18\\
10.5	15\\
11.2	15\\
11.2	-18\\
}--cycle;

\addplot[area legend, draw=none, fill=red, fill opacity=0.2, forget plot]
table[row sep=crcr] {%
x	y\\
13.075	-18\\
13.075	15\\
13.78	15\\
13.78	-18\\
}--cycle;
\end{axis}

\begin{axis}[%
width=0.952\fwidth,
height=0.410\fheight,
at={(0\fwidth,0\fheight)},
scale only axis,
xmin=0,
xmax=12,
xlabel style={font=\color{white!15!black}},
xlabel={$t$ [physical time]},
ymin=-18,
ymax=15,
ylabel style={font=\color{white!15!black}},
ylabel={$q(t),\ v(t)$},
axis background/.style={fill=white},
xmajorgrids,
ymajorgrids,
legend style={at={(0.55,0.03)}, anchor=south west, legend cell align=left, align=left, draw=white!15!black, nodes={scale=0.7, transform shape}},
legend columns=2, 
xlabel style={font={\scriptsize}},ylabel style={font=\scriptsize},ylabel shift={-0cm},ticklabel style={font=\scriptsize}
]
\addplot [color=mycolor1]
  table[row sep=crcr]{%
0	10\\
0.115	9.935131375\\
0.244999999999999	9.705577375\\
0.404999999999999	9.195457375\\
0.595000000000001	8.263507375\\
0.82	6.701878\\
1.085	4.225711375\\
1.385	0.591106375000004\\
1.42784312292855	-0.51519485476563\\
1.42784312292855	-2.97461164370512\\
1.46712357987742	0.487618019553551\\
1.78712357987742	3.8960934274499\\
2.06212357987742	6.02267135611082\\
2.29712357987742	7.25206460878469\\
2.48712357987742	7.84996338222314\\
2.63712357987742	8.07183372967456\\
2.75212357987742	8.0924544543873\\
2.87212357987742	7.97565073234843\\
3.01212357987742	7.66083772330308\\
3.18212357987742	7.02007128374802\\
3.38212357987742	5.90325841368323\\
3.62212357987742	4.04511496960549\\
3.89712357987742	1.22132164826641\\
3.99796074420031	-0.426808388213457\\
3.99796074420031	-2.67723403438744\\
4.03424715975103	0.405239408767844\\
4.33924715975103	3.30084520704264\\
4.59924715975103	5.04868597770312\\
4.81424715975103	5.99309272074928\\
4.98924715975103	6.42702963369384\\
5.12424715975103	6.55650671653678\\
5.23924715975103	6.52578325932892\\
5.35924715975103	6.35540299963375\\
5.50924715975103	5.9437751750148\\
5.68924715975103	5.15846478547206\\
5.90424715975103	3.80389802851822\\
6.15924715975103	1.60945460166601\\
6.31106660334263	-0.26162065800226\\
6.31106660334263	-2.40949804117373\\
6.35137073962436	0.403586034001144\\
6.63637073962436	2.80268674690364\\
6.88137073962436	4.22815732027595\\
7.08137073962435	4.95526176792682\\
7.24137073962435	5.25441732604752\\
7.36637073962435	5.31339198082932\\
7.48137073962435	5.23227066322857\\
7.61137073962435	4.98434405420164\\
7.77137073962435	4.45161161232234\\
7.96137073962435	3.49280933759066\\
8.19137073962435	1.85833095238916\\
8.39286187656759	-0.137579105561965\\
8.39286187656759	-2.16855668055263\\
8.433494319496	0.365317481242046\\
8.703494319496	2.3814428427385\\
8.933494319496	3.53480796549474\\
9.118494319496	4.08593331966823\\
9.263494319496	4.28319218510151\\
9.378494319496	4.29297937147964\\
9.493494319496	4.17302930785775\\
9.633494319496	3.85189464344851\\
9.803494319496	3.20345190809443\\
10.008494319496	2.04443438163803\\
10.248494319496	0.16368181407932\\
10.2664776224657	-0.336762870093057\\
10.2664776224657	-1.95163256923867\\
10.3056178993727	0.31621796267566\\
10.5656178993727	2.03528871706635\\
10.7806178993727	2.95590486973557\\
10.9506178993727	3.36280167068333\\
11.0806178993727	3.48266304787867\\
11.1956178993727	3.45049589116686\\
11.3156178993727	3.2786091624241\\
11.4656178993727	2.86509825149565\\
11.4856178993727	2.79328646337186\\
};
\addlegendentry{$q(t)$}

\addplot [color=mycolor2]
  table[row sep=crcr]{%
0	0\\
1.42784312292855	-13.9974352153923\\
1.42784312292855	12.6347411500197\\
1.53212357987742	11.5834356496761\\
3.99712357987742	-12.5982143503239\\
3.99796074420031	-9.24533174167783\\
3.99796074420031	11.3711372363155\\
4.14424715975103	9.91071450254031\\
6.31106660334263	-11.3338239242465\\
6.31106660334263	10.2341836161548\\
6.44637073962435	8.88387223825437\\
8.39286187656759	-10.1992134706011\\
8.39286187656759	9.21073686568865\\
8.548494319496	7.6633309685054\\
10.2664776224657	-9.18415695544899\\
10.2664776224657	8.28963534106107\\
10.4006178993727	6.95516059381035\\
11.4856178993727	-3.68868940618966\\
};
\addlegendentry{$v(t)$}

\end{axis}

\begin{axis}[%
width=1.228\fwidth,
height=1.228\fheight,
at={(-0.16\fwidth,-0.136\fheight)},
scale only axis,
xmin=0,
xmax=1,
ymin=0,
ymax=1,
axis line style={draw=none},
ticks=none,
axis x line*=bottom,
axis y line*=left,
legend style={legend cell align=left, align=left, draw=white!15!black},
xlabel style={font={\scriptsize}},ylabel style={font=\scriptsize},ylabel shift={-0cm},ticklabel style={font=\scriptsize}
]
\end{axis}
\end{tikzpicture}%

%% file: img/numerical_k5.tikz
\definecolor{mycolor1}{rgb}{0.00000,0.44700,0.74100}%
\definecolor{mycolor2}{rgb}{0.85000,0.32500,0.09800}%
\definecolor{mycolor3orange}{rgb}{0.9290, 0.6940, 0.1250}%
\definecolor{mycolor4bred}{rgb}{0.6350, 0.0780, 0.1840}%

\pgfplotsset{compat=1.13}
\setlength{\fwidth}{6.5cm}
\setlength{\fheight}{2.7cm}
%\pgfplotsset{every tick label/.append style={font=\tiny}} 
\pgfplotsset{
	every axis plot/.append style={line width=1.0pt},
	every axis plot post/.append style={
		every mark/.append style={line width=0.0pt}
		%every tick label/.append style={font=\tiny}
	}
}

\begin{tikzpicture}

\begin{axis}[%
width=0.951\fwidth,
height=\fheight,
at={(0\fwidth,0\fheight)},
scale only axis,
xmode=log,
xmin=10,
xmax=100000,
xminorticks=true,
xlabel style={font=\color{white!15!black}},
xlabel={Function Evaluations},
ymode=log,
ymin=0.0001,
ymax=100,
yminorticks=true,
ylabel style={font=\color{white!15!black}, nodes={scale=0.70, transform shape}},
ylabel={$E(1)$},
axis background/.style={fill=white},
xmajorgrids,
xminorgrids,
ymajorgrids,
yminorgrids,
legend style={legend cell align=left, align=left, draw=white!15!black,nodes={scale=0.9, transform shape}},
xlabel style={font={\scriptsize}},ylabel style={font=\scriptsize},ylabel shift={-0cm},ticklabel style={font=\scriptsize}
]
\addplot [color=mycolor3orange, mark=diamond, mark options={solid, black}]
  table[row sep=crcr]{%
10	24.5244580378635\\
50	2.79297092947111\\
100	2.55969454247095\\
500	0.828053432331639\\
1000	0.390910276824018\\
5000	0.0870792624907055\\
10000	0.0472950976119251\\
50000	0.0073907555535071\\
100000	0.00416025789417728\\
};
\addlegendentry{Explicit Euler}

\addplot [color=mycolor4bred, mark=o, mark options={solid, black}]
  table[row sep=crcr]{%
10	20.1209725005759\\
50	4.88402641655974\\
100	0.447879987749085\\
500	0.16414769779553\\
1000	0.0204564190605805\\
5000	0.0292442141407767\\
10000	0.0290110253671421\\
50000	0.000772170804680172\\
100000	0.00115225433363419\\
};
\addlegendentry{$\text{Runge}-\text{Kutta }4$}

\end{axis}

\begin{axis}[%
width=1.227\fwidth,
height=1.227\fheight,
at={(-0.16\fwidth,-0.135\fheight)},
scale only axis,
xmin=0,
xmax=1,
ymin=0,
ymax=1,
axis line style={draw=none},
ticks=none,
axis x line*=bottom,
axis y line*=left,
legend style={legend cell align=left, align=left, draw=white!15!black},
xlabel style={font={\scriptsize}},ylabel style={font=\scriptsize},ylabel shift={-0cm},ticklabel style={font=\scriptsize}
]
\end{axis}
\end{tikzpicture}%

%% file: img/sim3d.tikz
\pgfplotsset{compat=1.13}

\setlength{\fwidth}{6.5cm}
\setlength{\fheight}{2.7cm}
%\pgfplotsset{every tick label/.append style={font=\tiny}} 
\pgfplotsset{
	every axis plot/.append style={line width=1.0pt},
	every axis plot post/.append style={
		every mark/.append style={line width=0.0pt}
		%every tick label/.append style={font=\tiny}
	}
}

\usetikzlibrary{arrows}
\usetikzlibrary{shapes,snakes}
\usetikzlibrary{arrows.meta}
\begin{tikzpicture}

\begin{axis}[%
width=0.95\fwidth,
height=\fheight,
at={(-0.027\fwidth,-0\fheight)},
scale only axis,
unbounded coords=jump,
xmin=-0.586133543348177,
xmax=3.45541338359545,
tick align=outside,
xlabel style={font=\color{white!15!black}},
xlabel={$x$},
ymin=-0.685693807681318,
ymax=3.60967742714902,
ylabel style={font=\color{white!15!black}},
ylabel={$y$},
zmin=-0.141029874956004,
zmax=4.84628959900958,
zlabel style={font=\color{white!15!black}},
zlabel={$z$},
%view={-188.406055018886}{42.1153833168823},
view={-198}{43},
axis background/.style={fill=white},
axis x line*=bottom,
axis y line*=left,
axis z line*=left,
xmajorgrids,
ymajorgrids,
zmajorgrids,
legend style={at={(1.03,1)}, anchor=north west, legend cell align=left, align=left, draw=white!15!black},
xlabel style={font={\scriptsize}},ylabel style={font=\scriptsize},ylabel shift={-0cm},ticklabel style={font=\scriptsize}
]
%\draw [red,ultra thick,-{Latex[length=1.5mm, width=1.5mm]}] (1,1,4) --  (-0.1,-0.15,4);
\draw [->,red,ultra thick] (1,1,4) --    (0,-0.2670,4.0000);
\addplot3 [color=black, line width=1.2pt, draw=none, mark size=1.5pt, mark=o, mark options={solid, black}]
 table[row sep=crcr] {%
1	1	4\\
0.939216369375412	0.929085764271314	3.99798641588509\\
0.878432738750824	0.858171528542629	3.99194566353997\\
0.817649108126236	0.787257292813943	3.98187774296425\\
0.756865477501648	0.716343057085258	3.96778265415759\\
0.696081846877061	0.645428821356572	3.9496603971199\\
0.635298216252473	0.574514585627886	3.92751097185141\\
0.574514585627885	0.5036003498992	3.90133437835235\\
0.513730955003298	0.432686114170514	3.87113061662334\\
0.45294732437871	0.361771878441829	3.83689968666508\\
0.392163693754122	0.290857642713143	3.79864158847822\\
0.331380063129534	0.219943406984457	3.75635632206308\\
0.270596432504947	0.149029171255771	3.71004388741991\\
0.209812801880359	0.0781149355270854	3.65970428454865\\
0.149029171255771	0.00720069979839977	3.60533751344891\\
0.0906398689362478	0.0578761032866009	3.54931429096789\\
0.0298562383116598	0.125244565866262	3.48705181761333\\
0.055980363556662	0.223645565488371	3.38886329351496\\
0.113724842086582	0.291014028068032	3.31669141320505\\
0.171469320616502	0.358382490647694	3.24049236466411\\
0.229213799146422	0.425750953227355	3.16026614789177\\
0.286958277676342	0.493119415807016	3.07601276288779\\
0.344702756206262	0.560487878386677	2.98773220965181\\
0.402447234736182	0.627856340966338	2.8954244881836\\
0.460191713266102	0.695224803546	2.79908959848291\\
0.517936191796022	0.762593266125661	2.69872754055007\\
0.575680670325943	0.829961728705323	2.59433831438532\\
0.633425148855863	0.897330191284984	2.48592191998798\\
0.691169627385783	0.964698653864645	2.37347835735824\\
0.748914105915702	1.03206711644431	2.25700762649642\\
0.806658584445622	1.09943557902397	2.13650972740292\\
0.864403062975542	1.16680404160363	2.01198466007831\\
0.922147541505463	1.23417250418329	1.88343242452329\\
0.979892020035383	1.30154096676295	1.75085302073868\\
1.0376364985653	1.36890942934261	1.6142464487254\\
1.09538097709522	1.43627789192228	1.47361270848443\\
1.15312545562514	1.50364635450194	1.32895180001668\\
1.21086993415506	1.5710148170816	1.18026372332295\\
1.26861441268498	1.63838327966126	1.02754847840367\\
1.3263588912149	1.70575174224092	0.870806065258347\\
1.38410336974482	1.77312020482058	0.710036483887718\\
1.44184784827474	1.84048866740024	0.545239734291321\\
1.49959232680466	1.90785712997991	0.376415816468422\\
1.55733680533458	1.97522559255957	0.203564730418294\\
1.6150812838645	2.04259405513923	0.0266864761400125\\
1.67056196259091	2.10732142170013	0.137084009899427\\
1.72830644112083	2.17468988427979	0.302318226979606\\
1.78605091965075	2.24205834685945	0.463525275828355\\
1.84379539818067	2.30942680943911	0.62070515644483\\
1.90153987671059	2.37679527201877	0.773857868828094\\
1.95928435524051	2.44416373459844	0.922983412977219\\
2.01702883377043	2.5115321971781	1.0680817888921\\
2.07477331230035	2.57890065975776	1.20915299657183\\
};

\addplot3 [color=blue, dashed]
 table[row sep=crcr] {%
1	1	4\\
0.756865477501648	0.716343057085258	3.96778265415759\\
0.544122770315592	0.468143232034858	3.88673589351654\\
0.331380063129534	0.219943406984457	3.75635632206308\\
0.142859494113408	-0.0282414377376563	3.59959403342984\\
0.142859494113408	-1.07901919108169	3.59959403342984\\
0.142859494113408	-0.0094936180105325	3.59959403342984\\
-0.000535570619303094	0.158334028259741	3.45499548032541\\
-0.924795803136869	0.158334028259741	3.45499548032541\\
-0.00176417009400787	0.158334028259741	3.45499548032541\\
0.171469320616502	0.358382490647694	3.24049236466411\\
0.460191713266102	0.695224803546	2.79908959848291\\
0.748914105915702	1.03206711644431	2.25700762649642\\
1.06650873783026	1.40259366063245	1.54443297463332\\
1.38410336974482	1.77312020482058	0.710036483887718\\
1.6150812838645	2.04259405513923	0.0266864761400125\\
1.6236832009038	2.05262961076165	-0.241860716070098\\
1.6236832009038	2.05262961076165	-2.08539028197143\\
nan	nan	nan\\
1.6236832009038	2.05262961076165	-2.07573605931379\\
1.6236832009038	2.05262961076165	-0.0321115482461467\\
1.67056196259091	2.10732142170013	0.137084009899427\\
1.90153987671059	2.37679527201877	0.773857868828094\\
2.10364555156531	2.61258489104759	1.27817841232324\\
};

\addplot3[area legend, line width=0.5pt, draw=black, fill=black, fill opacity=0.15]
table[row sep=crcr] {%
x	y	z\\
0	0	0\\
3.5	0	0\\
3.5	3.5	0\\
0	3.5	0\\
}--cycle;

%\addplot3[area legend, line width=0.5pt, draw=black, fill=black, fill opacity=0.15]
%table[row sep=crcr] {%
%x	y	z\\
%0	0	0\\
%3.5	0	0\\
%3.5	3.5	0\\
%3.5	3.5	0\\
%}--cycle;

\addplot3[area legend, line width=0.5pt, draw=black, fill=black, fill opacity=0.15]
table[row sep=crcr] {%
x	y	z\\
0	0	0\\
3.5	0	0\\
3.5	0	5\\
0	0	5\\
}--cycle;

\addplot3[area legend, line width=0.5pt, draw=black, fill=black, fill opacity=0.15]
table[row sep=crcr] {%
x	y	z\\
0	0	0\\
0	3.5	0\\
0	3.5	5\\
0	0	5\\
}--cycle;

\end{axis}

\begin{axis}[%
width=1.227\fwidth,
height=1.227\fheight,
at={(-0.16\fwidth,-0.135\fheight)},
scale only axis,
xmin=0,
xmax=1,
ymin=0,
ymax=1,
axis line style={draw=none},
ticks=none,
axis x line*=bottom,
axis y line*=left,
legend style={legend cell align=left, align=left, draw=white!15!black},
xlabel style={font={\scriptsize}},ylabel style={font=\scriptsize},ylabel shift={-0cm},ticklabel style={font=\scriptsize}
]
\end{axis}
\end{tikzpicture}%

%% file: img/ocp3d.tikz
\pgfplotsset{compat=1.13}

\setlength{\fwidth}{6.5cm}
\setlength{\fheight}{3.2cm}
%\pgfplotsset{every tick label/.append style={font=\tiny}} 
\pgfplotsset{
	every axis plot/.append style={line width=1.0pt},
	every axis plot post/.append style={
		every mark/.append style={line width=0.0pt}
		%every tick label/.append style={font=\tiny}
	}
}

%\usetikzlibrary{arrows}
%\usetikzlibrary{shapes,snakes}
%\usetikzlibrary{arrows.meta}
%\usetikzlibrary{arrows.meta,arrows}

\begin{tikzpicture}
\begin{axis}[%
width=0.95\fwidth,
height=\fheight,
at={(-0.034\fwidth,0\fheight)},
scale only axis,
xmin=-1.30005544113432,
xmax=6.44787844596894,
tick align=outside,
xlabel style={font=\color{white!15!black}},
xlabel={$x$},
ymin=-2.25877531951624,
ymax=5.48915857078767,
ylabel style={font=\color{white!15!black}},
ylabel={$y$},
zmin=-0.248340356264137,
zmax=4.68216302465247,
zlabel style={font=\color{white!15!black}},
zlabel={$z$},
%view={162.300000130148}{22.2000003983238},
%view={-195.300000130148}{53.2000003983238},
view={-188.0}{45},
%view={0}{90},
axis background/.style={fill=white},
axis x line*=bottom,
axis y line*=left,
axis z line*=left,
xmajorgrids,
ymajorgrids,
zmajorgrids,
legend style={at={(1.03,1)}, anchor=north west, legend cell align=left, align=left, draw=white!15!black},
xlabel style={font={\scriptsize}},ylabel style={font=\scriptsize},ylabel shift={-0cm},ticklabel style={font=\scriptsize}
]
%\usetikzlibrary{arrows.meta,arrows}

\addplot3 [color=black, line width=1.0pt, draw=none, mark=o, mark options={solid, black}]
 table[row sep=crcr] {%
1	1	4\\
0.967978244220108	0.964878666414156	3.99999881185458\\
%0.934910626882646	0.929319568209113	3.99999637940544\\
0.900797147988139	0.892741957610885	3.99999264129766\\
%0.865637807537189	0.855132467189071	3.99998753022701\\
0.829432605530495	0.816486526896694	3.99998097203457\\
%0.792181541968864	0.776801829481409	3.99997288460786\\
0.753884616853241	0.736076983414724	3.99996317653384\\
%0.714541830184727	0.694311058036535	3.99995174542977\\
0.674153181964623	0.651503387158293	3.99993847584915\\
%0.63271867219447	0.607653470379704	3.99992323661868\\
0.59023830087611	0.562760918071408	3.99990587739864\\
%0.546712068011762	0.516825418277068	3.99988622416156\\
0.502139973604137	0.469846715607436	3.99986407312754\\
%0.456522017656586	0.421824597146342	3.99983918243625\\
0.409858200173334	0.372758882679951	3.99981126039123\\
%0.362148521159814	0.322649417709718	3.99977994831057\\
0.313392980623231	0.271496068323587	3.99974479448995\\
%0.263591578573517	0.219298717346263	3.99970521266319\\
0.212744315025183	0.166057261393301	3.99966041140466\\
%0.160851190001542	0.111771608578789	3.9996092636775\\
0.107912203548719	0.0564416767061608	3.99955003602973\\
%0.0539273561348166	6.7391837074382e-05	3.99947971919999\\
0.0497278201708884	0.0359636428626544	3.99939180489156\\
%0.100490111688682	0.0889507331607016	3.99928355994351\\
0.152295517562888	0.142980490787591	3.99915138264026\\
%0.20514394001684	0.198052791724162	3.99899030476635\\
0.259035275971591	0.254167504554651	3.99879304388713\\
%0.313969416653289	0.311324489862032	3.99854795894162\\
0.369946247165708	0.369523599561338	3.99823293896543\\
%0.426965646022418	0.42876467616205	3.99779102018342\\
0.485027484633589	0.489047551950017	3.99689179786715\\
0.544131626741956	0.550372048077995	3.99410363673608\\
%0.604277927801756	0.612737973552214	3.98931915982303\\
0.665466234293554	0.676145124100249	3.9825059510787\\
%0.727696382966854	0.74059328090303	3.97364849404213\\
0.790968200001179	0.806082209170805	3.96273772029286\\
%0.855281500074837	0.872611656539269	3.94976769831119\\
0.920636085328939	0.940181351257737	3.93473425788296\\
%0.987031744212189	1.00879100013587	3.91763431670888\\
1.05446825018962	1.07844028620899	3.8984655120504\\
%1.12294536029557	1.14912886607412	3.87722598235783\\
1.19246281350781	1.22085636683868	3.85391422966705\\
%1.26302032891565	1.29362238261198	3.82852902864235\\
1.33461760364982	1.36742647045365	3.8010693640986\\
%1.40725431053601	1.44226814567408	3.77153438671852\\
1.48093009542659	1.51814687635706	3.73992338084387\\
%1.55564457415596	1.5950620769431	3.70623574054152\\
1.63139732905401	1.67301310067083	3.67047095150252\\
%1.70818790493827	1.75199923062085	3.63262857715664\\
1.78601580448833	1.83201966903537	3.59270824790341\\
%1.86488048288442	1.91307352449402	3.55070965269686\\
1.94478134156482	1.99515979639965	3.50663253244466\\
%2.02571772092191	2.07827735605605	3.46047667483565\\
2.10768889171234	2.1624249233821	3.41224191031614\\
%2.19069404489861	2.24760103797302	3.36192810901276\\
2.27473227956406	2.33380402274395	3.30953517845567\\
2.35980258844264	2.42103193769869	3.25506306199942\\
%2.44590384047373	2.50928252034425	3.19851173787896\\
2.53303475961062	2.59855310771065	3.13988121884162\\
%2.62119389886405	2.68884053251784	3.07917155235899\\
2.71037960822438	2.780140982142	3.01638282142098\\
%2.80058999462246	2.87244980258616	2.95151514594573\\
2.89182287140064	2.96576121852926	2.88456868486628\\
%2.98407569375781	3.0600679204037	2.81554363898723\\
3.07734547512672	3.15536043099409	2.74444025474503\\
%3.17162867713716	3.25162608544045	2.67125882905826\\
3.26692106218775	3.34884728388734	2.5959997155247\\
3.36321749174318	3.44699824490994	2.51866333232011\\
%3.46051164349853	3.5460382691206	2.43925017229275\\
3.55879560295635	3.64589542650755	2.35776081594971\\
%3.65805925230607	3.74641798297491	2.2741959483331\\
3.75828931518135	3.84721197782025	2.1885563812464\\
%3.85946778011634	3.94757650608117	2.10084308302062\\
3.96156911504393	4.04725851402815	2.01105721919905\\
%4.06455490786833	4.14618108573501	1.91920020953408\\
4.16836242952043	4.24429884472105	1.82527381025564\\
%4.27287731079464	4.34158718600465	1.72928023722588\\
4.37786380878679	4.43803475451697	1.63122235882635\\
4.48283331690251	4.53364003029467	1.5311040158835\\
%4.58722327145072	4.62841087166977	1.42893059359203\\
4.6910027989686	4.72236684919965	1.32471015402977\\
%4.79431828160147	4.81554641846665	1.21845603828446\\
4.89727316973706	4.90802751947218	1.11019446257941\\
5.0000000039834	5.00000000852747	1.00000003421128\\
};

\addplot3[color=blue, dashed,line width = 0.5pt]
 table[row sep=crcr] {%
1	1	4\\
0.0539273561348166	6.7391837074382e-05	3.99947971919999\\
-0.981634982817951	-0.0562542298823292	3.99947971922479\\
-1.35438014561707	-0.650398316693271	3.99947971923251\\
-1.52723786462239	-1.06829782720609	3.99947971923884\\
-1.5694896102719	-1.35789088279626	3.99947971924445\\
-1.50555084290819	-1.55133497063981	3.99947971925007\\
-1.31453644059715	-1.64015380872809	3.99947971925639\\
-0.976874458376686	-1.5738624766606	3.99947971926342\\
-0.497579995565284	-1.31360951531585	3.99947971927114\\
-0.0152870082238268	-0.937221639690498	3.99947971927817\\
4.30548670316e-06	-0.665098461201984	3.99947971928635\\
8.73586098837364e-06	-0.0159806629276424	3.99947971930427\\
0.727696382966854	0.74059328090303	3.97364849404213\\
1.12294536029557	1.14912886607412	3.87722598235783\\
1.63139732905401	1.67301310067083	3.67047095150252\\
2.27473227956406	2.33380402274395	3.30953517845567\\
3.07734547512672	3.15536043099409	2.74444025474503\\
4.27287731079464	4.34158718600465	1.72928023722588\\
5.0000000039834	5.00000000852747	1.00000003421128\\
};

%\draw [->,red,ultra thick] (1,1,4) --    (0,-0.42,4.0000);
\draw [->,red,ultra thick] (1,1,4) --    (-0.2,-0.95,4.0000);
\addplot3[area legend,line width=0.5pt, draw=black, fill=black, fill opacity=0.15]
table[row sep=crcr] {%
x	y	z\\
0	0	0\\
6	0	0\\
6	6	0\\
0	6	0\\
}--cycle;

%
%\addplot3[area legend, draw=black, fill=black, fill opacity=0.2]
%table[row sep=crcr] {%
%x	y	z\\
%0	0	0\\
%5	0	0\\
%5	5	0\\
%0	5	0\\
%}--cycle;

\addplot3[area legend,line width=0.5pt, draw=black, fill=black, fill opacity=0.15]
table[row sep=crcr] {%
x	y	z\\
0	0	0\\
6	0	0\\
6	0	5\\
0	0	5\\
}--cycle;

\addplot3[area legend,line width=0.5pt, draw=black, fill=black, fill opacity=0.15]
table[row sep=crcr] {%
x	y	z\\
0	0	0\\
0	6	0\\
0	6	5\\
0	0	5\\
}--cycle;

\end{axis}

\begin{axis}[%
width=1.227\fwidth,
height=1.227\fheight,
at={(-0.16\fwidth,-0.135\fheight)},
scale only axis,
xmin=0,
xmax=1,
ymin=0,
ymax=1,
axis line style={draw=none},
ticks=none,
axis x line*=bottom,
axis y line*=left,
legend style={legend cell align=left, align=left, draw=white!15!black},
xlabel style={font={\scriptsize}},ylabel style={font=\scriptsize},ylabel shift={-0cm},ticklabel style={font=\scriptsize}
]
\end{axis}
\end{tikzpicture}%